\pdfoutput=1
\documentclass[hidelinks,onefignum,onetabnum]{siamart251216}

\usepackage{pgfplots}
\pgfplotsset{compat=1.18}
\usepackage{bm}
\usepackage{booktabs}
\usepgfplotslibrary{groupplots}
\usepackage{amsfonts,amssymb}
\usepackage{comment}
\usepackage{graphicx}
\usepackage{epstopdf}
\usepackage[T1]{fontenc}
\usepackage[utf8]{inputenc} 
\usepackage{tikz}
\usepackage{orcidlink}
\ifpdf
  \DeclareGraphicsExtensions{.eps,.pdf,.png,.jpg}
\else
  \DeclareGraphicsExtensions{.eps}
\fi

\newsiamremark{remark}{Remark}
\newsiamremark{hypothesis}{Hypothesis}
\crefname{hypothesis}{Hypothesis}{Hypotheses}
\newsiamthm{claim}{Claim}
\newsiamremark{fact}{Fact}
\crefname{fact}{Fact}{Facts}

\headers{Fast matrix--vector products in tensor-train format}{A.\ Bucci, S.\ Güttel, J.\,W.\ Pearson}

\title{Fast matrix--vector products in tensor-train format via structured sketching\thanks{
\funding{This work was funded by the Engineering and Physical Sciences Research Council (EPSRC) UK grant EP/Z533786/1.}}}

\author{Alberto Bucci\thanks{School of Mathematics and Maxwell Institute for Mathematical Sciences, The University of Edinburgh, James Clerk Maxwell Building, The King’s Buildings, Edinburgh, EH9 3FD, United
Kingdom  
  (\email{abucci2@ed.ac.uk} \orcidlink{0009-0002-7550-6025}, \email{j.pearson@ed.ac.uk} \orcidlink{0000-0002-6063-1766})}
\and Stefan Güttel\thanks{Department of Mathematics, The University of Manchester, Oxford Road, Manchester, M13 9PL, United Kingdom (\email{stefan.guettel@manchester.ac.uk} \orcidlink{0000-0003-1494-4478})}
\and John W.~Pearson\footnotemark[2]
}

\usepackage{amsopn}
\DeclareMathOperator{\diag}{diag}

\ifpdf
\hypersetup{
  pdftitle={Fast matrix–vector products in tensor-train format via structured sketching},
  pdfauthor={A. Bucci, S. Güttel, and J. W. Pearson}
}
\fi

\begin{document}

\maketitle

\begin{abstract}
We propose a randomized, streamable algorithm that computes a tensor-train approximation of a matrix--vector product without ever forming the full intermediate tensor. We develop three contraction strategies for streaming tensor-train approximation (STTA): sequential, fused, and pre-sketched. We focus on the regime of large operator mode sizes and show that pre-sketching can provide significant speed-ups when a lower-rank approximation of the input is sought, by pre-computing all sketch-related contractions with the operator. The pre-sketched strategy can also accelerate successive randomized compression (SRC), although to a lesser extent since its second sweep depends on the input vector and cannot be pre-computed. We compare the three strategies in terms of computational cost and memory, and demonstrate their performance on synthetic and application-driven examples where rank reduction is relevant.
\end{abstract}



\section{Introduction}
\label{sec:intro}

We study randomized, streamable techniques for constructing low-rank approximations of the
matrix--vector product (matvec) in the tensor-train (TT) format between a
TT-operator $A$ and a TT-vector $v$, often referred to as the MPO--MPS
application in the physics literature
\cite{oseledets2011tensor, schollwock2011density}. The TT matvec is a fundamental
building block in many modern tensor-based computational workflows, including
quantum chemistry \cite{szalay2015tensor}, quantum many-body simulation
\cite{khoromskij2011d}, large-scale discretizations of partial differential
equations \cite{bachmayr2016tensor}, and high-dimensional model reduction
\cite{Nouy2017}. Accelerating the matvec routine is therefore of critical importance.

The difficulty is well known. If $A$ has TT ranks $r^A$ and $v$ has TT ranks
$r^v$, the cores of $Av$ carry TT ranks $r^Ar^v$. In almost every application such rank increase is undesirable: one instead wants an approximation of $Av$ of 
much smaller rank $s$, either because the result is genuinely of low rank, or
because a rank cap must be enforced to keep a recurrence affordable. Forming the
product and then compressing it would incur a large amount of work, with the compression step being  more
expensive than the formation, leading to severe computational and memory bottlenecks.

In this work, we demonstrate that \emph{streaming tensor-train approximation} (STTA)
\cite{kressner2023streaming}, an extension of the generalized Nystr\"om method 
\cite{clarkson2009numerical, nakatsukasa2020fast,tropp2017practical} to the TT format, can be utilized to 
compute a low-rank approximation of $Av$ without ever forming the full product. 
While this capability might be anticipated given recent parallel advances like the 
successive randomized compression (SRC) \cite{camano2026successive}, a tensor adaptation of the randomized singular value decomposition (SVD) \cite{halko2011finding} that is closely related to randomized TT-rounding algorithms \cite{aldaas2023randomized}, our contribution is particularly tailored to regimes of moderate tensor order~$d$ and large mode sizes. In this 
setting, the cost is driven by the per-mode dimensions rather than by the order, 
requiring specific algorithmic designs to achieve high efficiency. 

To this end, we characterize the optimal ordering of tensor contractions required to obtain 
the STTA approximation without instantiating the full product. We propose and analyze 
three distinct approaches: \emph{sequential}, \emph{fused}, and \emph{pre-sketched}. 
They all compute the same object but differ in the order in which the core of the operator, 
the core of the vector, and the sketching matrices are contracted at each site.
The sequential algorithm contracts the operator core with the vector core
first, producing the cores of $Av$ one at a time and compressing each of them
before moving on. The fused algorithm contracts the vector core into the
running sketches instead, avoiding the formation of any core of $Av$. 
The pre-sketched algorithm does the same, but contracts the sketch into the cores of the operator beforehand. The preferable strategy depends on the dimensions and TT ranks involved, as well as on whether computational cost or memory usage is the primary concern. In the regime considered here, however, the main factor is the relative size of the target rank~$s$ and the TT ranks~$r^v$ of the vector. Here and below, comparisons between a scalar and the tuple \(r^v\) are understood componentwise.

When $s\lesssim r^v$, that is, when one seeks an approximation whose rank is smaller 
than that of the input vector, the pre-sketched algorithm emerges as the most 
economical choice, and it represents the main algorithmic novelty of this work.  
Because STTA relies on sketching operators that are fixed, random, and drawn 
independently of the data, the sketch can be contracted into the cores of $A$ 
\emph{once and for all} in both the forward and backward sweeps. The step is an instance of the
associativity $(Y^{\top}A)v=Y^{\top}(Av)$. Subsequent contractions then work with an operator of which the row modes have already been replaced by the much smaller sketch dimension~$s$. Since this compressed operator is independent of $v$, it can be reused when the same operator is applied repeatedly, as in tensor-train adaptations of the power method \cite{oseledets2011tensor}, GMRES~\cite{TT-GMRES}, LSQR~\cite{paige1982lsqr}, and LSMR~\cite{fong2011lsmr}, as well as their randomized counterparts~\cite{bucci2025randomized,nakatsukasa2024fast}, and more generally in iterative Krylov methods.

\paragraph{Related work} Our approach builds on a line of recent work on
randomized compression in the TT format \cite{aldaas2023randomized,aldaas2026adaptive,camano2026successive,kressner2023streaming,meng2026recursive}.
The closest relative is the successive randomized compression approach of
\cite{camano2026successive}, which also avoids forming the full product. The
authors of \cite{camano2026successive} primarily target the regime of large $d$
and small mode sizes typical of quantum many-body physics. Both STTA and SRC perform two sweeps over the the cores of $A$, though they
differ crucially in the nature of the projections. STTA uses fixed oblique
projections for both sweeps. SRC, on the other hand, contracts the operator core during the
second sweep against an
orthonormal factor obtained from a data-dependent QR factorization.
In the large-mode regime considered here, the straightforward implementation of SRC described in \cite{camano2026successive} can be more expensive than the proposed STTA-based approach, depending on the relative sizes of~$s$ and~$r^v$. Nevertheless, when $s$ is sufficiently small compared to $r^v$, adapting the pre-sketching strategy to its first sweep and using a suitable contraction ordering in the second sweep brings its computational cost to a level comparable to that of the proposed STTA-based approach.
The trade-off is that SRC's orthogonal projections yield
slightly higher accuracy than the oblique projections of STTA; however,
available error analyses and empirical observations
\cite{aldaas2023randomized, cazeaux2026linear, kressner2023streaming} show that
this accuracy gap is prominent mainly for large $d$, whereas in our target
regime (small $d$, large mode sizes) the difference is often marginal.

Deterministic alternatives are widely used in practice and also avoid explicitly forming the product. The \emph{zip-up} method
\cite{paeckel2019time, stoudenmire2010minimally} sweeps once through the network, contracting the cores of $A$ and $v$ with the incoming carry at each site and truncating to the target rank before moving on, so that no core of
$Av$ is ever held at the full intermediate rank $r^A r^v$. Alternating
optimization schemes, such as DMRG-type matrix--vector multiplication and
AMEn \cite{dolgov2014alternating, oseledets2011dmrg}, instead fit a rank-$s$
tensor to $Av$ iteratively. Both approaches are often very accurate in
practice. However, neither provides an a priori guarantee: in the zip-up
method, the truncation at one site is decided before the environment to its
right has been formed, while the fitting schemes converge to a stationary
point of a nonconvex problem. Consequently, neither approach provides an
a priori bound on the error relative to the best rank-$s$ approximation of
$Av$, of the type available for the sketching-based methods considered here.

The approaches just discussed share a further limitation that is central to our setting:
they are not streamable. In particular, zip-up, alternating schemes,
and SRC cannot process a sum of TT tensors term by term while discarding
each term after its contribution has been incorporated. In contrast, the sketches used by STTA make the method naturally streamable: each term can be sketched as it becomes available and then discarded, without ever forming the full product, and the final approximation is reconstructed from the accumulated sketches.

Structured sketching operators for the TT format comprise an active subject area in their
own right. We use Khatri--Rao maps throughout, although other sketching maps, such as TT-Gaussian operators and intermediate TT-stack constructions, can also be employed \cite{ahle2020oblivious,camano2025faster,cazeaux2026linear}.
\paragraph{Contributions} The main contributions of this work are the following:
\begin{itemize}
  \item We show that STTA can be seamlessly adapted to compute a low-rank TT approximation of 
  $Av$ without ever forming and storing the full product. We propose and analyze three 
  contraction strategies to achieve this: sequential, fused, and pre-sketched.
  \item We characterize the computational and memory trade-offs of these three routes. 
  The choice in operations is settled by comparing $s$ with $r^v$: the pre-sketched route 
  is cheaper below that threshold, while the sequential route is better above the threshold.
  In memory, routes that avoid forming $Av$ (fused and pre-sketched) are vastly superior over 
  a much wider range of target ranks (\cref{sec:memory}). 
  \item We identify the pre-sketched algorithm as particularly effective in the regime $s \ll r^v$. By absorbing the sketch into the operator across both sweeps, the weight with which the cores of $A$ are visited drops from $(r^v)^2$ (sequential) to~$s^2$. We show that the resulting sketched operator can be recycled in iterative applications, and that in the same regime SRC can also be accelerated by exploiting the pre-sketching idea, achieving comparable computational costs.
 \item We validate our claims through synthetic benchmarks and two application scenarios in which $s \ll r^v$ arises naturally from the problem structure.
First, in a \emph{near-resonant solve}, truncating $v$ before applying $A$ destroys
accuracy by discarding components that are amplified by the operator. Our method
instead safely targets a lower-rank approximation without forming the full product.
Second, in a \emph{streaming uncertainty propagation} task, we accumulate an
ensemble of highly noisy realizations. Being streamable, our method can discard
each realization after sketching, avoiding a storage that grows with the length of the stream.
\end{itemize}

\paragraph{Outline} \Cref{sec:preliminaries} fixes the notation, recalls the TT
format, and summarizes STTA. \Cref{sec:fastmatvec} presents the three algorithms, then analyzes
the operation count, the memory footprint, and the range of validity, comparing these
with SRC. \Cref{sec:experiments} reports the
numerical experiments, and \cref{sec:conclusions} summarizes our findings.

\section{Notation and preliminaries}
\label{sec:preliminaries}

We briefly recall the tensor-train  format and introduce the randomized sketching tools used throughout the paper, assuming a basic familiarity with tensor-train arithmetic. We adopt standard TT notation and refer to the physics literature for the equivalent matrix product operator (MPO) and matrix product state (MPS) terminology.

\subsection{Tensor-train operators and vectors}

Let 
$A \in \mathbb{R}^{m_1 \times \dots \times m_d \times n_1 \times \dots \times n_d}$
be a tensor-train operator of order $d$, where $m_1,\dots,m_d$ and
$n_1,\dots,n_d$ denote the row and column mode sizes, respectively.
The operator is represented in TT format by a sequence of cores with
TT ranks $r^A = (r_0^A, r_1^A, \dots, r_{d-1}^A, r_d^A)$, where $r_0^A=r_d^A=1$. 
For convenience, we introduce the shorthand
\[
M := m_1 \cdots m_d,
\qquad
N := n_1 \cdots n_d,
\]
so that $A$ may be interpreted as a matrix in $\mathbb{R}^{M \times N}$,
while all computations are carried out in compressed TT form.

\begin{figure}[ht]
    \centering
\begin{tikzpicture}[baseline=-0.5ex, scale=0.60]

\tikzset{
  tensor/.style={circle, draw, thick, minimum size=4mm},
  green/.style={tensor, fill=green!60},
  blue/.style={tensor, fill=cyan!60},
  yellow/.style={tensor, fill=orange!80},
  wire/.style={thick},
}

\def\n{7}
\def\m{5}
\pgfmathtruncatemacro{\nplus}{\n+1}
\pgfmathtruncatemacro{\nminus}{\n-1}
\pgfmathtruncatemacro{\nplusm}{\n+\m}
\pgfmathtruncatemacro{\nplusmminus}{\n+\m-1}
\pgfmathtruncatemacro{\nplusmhalf}{\n+(\m-1)/2}
\pgfmathtruncatemacro{\nplusmhalfminus}{\n-1+(\m-1)/2}
\pgfmathtruncatemacro{\nhalfminus}{(\n)/2-1}
\pgfmathtruncatemacro{\nplusmhalfplus}{\n+(\m+3)/2}
\pgfmathtruncatemacro{\nhalf}{(\n-1)/2}
\pgfmathtruncatemacro{\nhalfplus}{(\n+3)/2}

\foreach \i in {1,...,\nhalf} {
  \node[blue] (b\i) at (\i,0) {};
}

\foreach \i in {\nhalfplus,...,\n} {
  \node[blue] (b\i) at (\i,0) {};
}
\node at (\nhalf+1,0) {$\cdots$};

\foreach \i in {1,...,\nhalfminus} {
  \draw[wire] (b\i) -- (b\the\numexpr\i+1\relax);
}

\foreach \i in {\nhalfplus,...,\nminus} {
  \draw[wire] (b\i) -- (b\the\numexpr\i+1\relax);
}
\draw[wire] (\nhalf+0.35,0) -- (\nhalf+0.65,0);
\draw[wire] (\nhalf+1.35,0) -- (\nhalf+1.65,0);

\foreach \i in {1,..., \nhalf} {
    \draw [wire] (b\i) -- (\i, 0.8);
}

\foreach \i in {\nhalfplus,..., \n} {
    \draw [wire] (b\i) -- (\i, 0.8);
}

\foreach \i in {1,..., \nhalf} {
    \draw [wire] (b\i) -- (\i, -0.8);
}

\foreach \i in {\nhalfplus,..., \n} {
    \draw [wire] (b\i) -- (\i, -0.8);
}
\node at (0,0) {$A=$};
\end{tikzpicture}

    \caption{Tensor diagram representation of $A$.}
    \label{fig:tt-operator}
\end{figure}
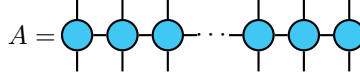

The TT representation of $A$ is given by a sequence of cores
\[
G_k \in \mathbb{R}^{r^A_{k-1} \times m_k \times n_k \times r^A_k},
\qquad k=1,\dots,d,
\]
such that
\[
[A]_{i_1,\dots,i_d,j_1,\dots,j_d}
=
\sum_{\alpha_1=1}^{r^A_1}\cdots
\sum_{\alpha_{d-1}=1}^{r^A_{d-1}}
G_1(i_1,j_1,\alpha_1)
G_2(\alpha_1,i_2,j_2,\alpha_2)
\cdots
G_d(\alpha_{d-1},i_d,j_d).
\]
This structure is illustrated in Figure~\ref{fig:tt-operator}.

\begin{figure}[ht]
\centering
   \begin{tikzpicture}[baseline=-0.5ex, scale=0.6]

\tikzset{
  tensor/.style={circle, draw, thick, minimum size=4mm},
  green/.style={tensor, fill=green!60},
  blue/.style={tensor, fill=cyan!60},
  yellow/.style={tensor, fill=orange!80},
  wire/.style={thick},
}

\def\n{7}
\def\m{5}
\pgfmathtruncatemacro{\nplus}{\n+1}
\pgfmathtruncatemacro{\nminus}{\n-1}
\pgfmathtruncatemacro{\nplusm}{\n+\m}
\pgfmathtruncatemacro{\nplusmminus}{\n+\m-1}
\pgfmathtruncatemacro{\nplusmhalf}{\n+(\m-1)/2}
\pgfmathtruncatemacro{\nplusmhalfminus}{\n-1+(\m-1)/2}
\pgfmathtruncatemacro{\nhalfminus}{(\n)/2-1}
\pgfmathtruncatemacro{\nplusmhalfplus}{\n+(\m+3)/2}
\pgfmathtruncatemacro{\nhalf}{(\n-1)/2}
\pgfmathtruncatemacro{\nhalfplus}{(\n+3)/2}

\foreach \i in {1,...,\nhalf} {
  \node[yellow] (y\i) at (\i,0) {};
}

\foreach \i in {\nhalfplus,...,\n} {
  \node[yellow] (y\i) at (\i,0) {};
}
\node at (\nhalf+1,0) {$\cdots$};

\foreach \i in {1,...,\nhalfminus} {
  \draw[wire] (y\i) -- (y\the\numexpr\i+1\relax);
}

\foreach \i in {\nhalfplus,...,\nminus} {
  \draw[wire] (y\i) -- (y\the\numexpr\i+1\relax);
}
\draw[wire] (\nhalf+0.35,0) -- (\nhalf+0.65,0);
\draw[wire] (\nhalf+1.35,0) -- (\nhalf+1.65,0);

\foreach \i in {1,..., \nhalf} {
    \draw [wire] (y\i) -- (\i, 0.8);
}

\foreach \i in {\nhalfplus,..., \n} {
    \draw [wire] (y\i) -- (\i, 0.8);
}

\node at (0,0) {$v=$};

\end{tikzpicture}
    \caption{Tensor diagram representation of $v$.}
    \label{fig:tt-vector}
\end{figure}
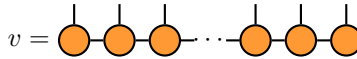

Let $v \in \mathbb{R}^{n_1 \times \dots \times n_d}$ be a TT-vector with TT-ranks $r^v = (r_0^v, r_1^v, \dots, r_{d-1}^v, r_d^v)$. The TT representation of the vector $v$ is defined by cores
\[
V_k \in \mathbb{R}^{r^v_{k-1} \times n_k \times r^v_k},
\qquad k=1,\dots,d,
\]
with $r^v_0=r^v_d=1$, such that
\[
[v]_{j_1,\dots,j_d}
=
\sum_{\beta_1=1}^{r^v_1}\cdots
\sum_{\beta_{d-1}=1}^{r^v_{d-1}}
V_1(j_1,\beta_1)
V_2(\beta_1,j_2,\beta_2)
\cdots
V_d(\beta_{d-1},j_d).
\]
A schematic representation is shown in Figure~\ref{fig:tt-vector}.

\subsection{TT matrix--vector product}

\begin{figure}[ht]
\centering
   \begin{tikzpicture}[baseline=-0.5ex, scale=0.48]

\tikzset{
  tensor/.style={circle, draw, thick, minimum size=3mm},
  green/.style={tensor, fill=green!60},
  blue/.style={tensor, fill=cyan!60},
  yellow/.style={tensor, fill=orange!80},
  wire/.style={thick},
}

\def\n{7}
\def\m{5}
\pgfmathtruncatemacro{\nplus}{\n+1}
\pgfmathtruncatemacro{\nminus}{\n-1}
\pgfmathtruncatemacro{\nplusm}{\n+\m}
\pgfmathtruncatemacro{\nplusmminus}{\n+\m-1}
\pgfmathtruncatemacro{\nplusmhalf}{\n+(\m-1)/2}
\pgfmathtruncatemacro{\nplusmhalfminus}{\n-1+(\m-1)/2}
\pgfmathtruncatemacro{\nhalfminus}{(\n)/2-1}
\pgfmathtruncatemacro{\nplusmhalfplus}{\n+(\m+3)/2}
\pgfmathtruncatemacro{\nhalf}{(\n-1)/2}
\pgfmathtruncatemacro{\nhalfplus}{(\n+3)/2}

\foreach \i in {1,...,\nhalf} {
  \node[blue] (b\i) at (\i,1) {};
}

\foreach \i in {\nhalfplus,...,\n} {
  \node[blue] (b\i) at (\i,1) {};
}

\foreach \i in {1,...,\nhalf} {
  \node[yellow] (y\i) at (\i,0) {};
}

\foreach \i in {\nhalfplus,...,\n} {
  \node[yellow] (y\i) at (\i,0) {};
}

\node at (\nhalf+1,0) {$\cdots$};
\node at (\nhalf+1,1) {$\cdots$};

\foreach \i in {1,...,\nhalfminus} {
  \draw[wire] (y\i) -- (y\the\numexpr\i+1\relax);
}

\foreach \i in {\nhalfplus,...,\nminus} {
  \draw[wire] (y\i) -- (y\the\numexpr\i+1\relax);
}

\foreach \i in {1,...,\nhalfminus} {
  \draw[wire] (b\i) -- (b\the\numexpr\i+1\relax);
}

\foreach \i in {\nhalfplus,...,\nminus} {
  \draw[wire] (b\i) -- (b\the\numexpr\i+1\relax);
}

\draw[wire] (\nhalf+0.35,0) -- (\nhalf+0.65,0);
\draw[wire] (\nhalf+1.35,0) -- (\nhalf+1.65,0);

\draw[wire] (\nhalf+0.35,1) -- (\nhalf+0.65,1);
\draw[wire] (\nhalf+1.35,1) -- (\nhalf+1.65,1);

\foreach \i in {1,..., \nhalf} {
    \draw [wire] (y\i) -- (b\i);
}

\foreach \i in {\nhalfplus,..., \n} {
    \draw [wire] (y\i) -- (b\i);
}

\foreach \i in {1,..., \nhalf} {
    \draw [wire] (\i,1.7) -- (b\i);
}

\foreach \i in {\nhalfplus,..., \n} {
    \draw [wire] (\i,1.7) -- (b\i);
}

\foreach \i in {\nhalfplus,..., \n} {
    \draw [wire] (y\i) -- (b\i);
}

\node at (-0.3,0.5) {$Av=$};

\end{tikzpicture}
    \caption{Tensor diagram representation of the $Av$ product.}
    \label{fig:tt-Av}
\end{figure}
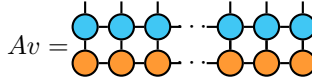

The TT matrix--vector product $Av$, depicted in Figure~\ref{fig:tt-Av}, yields a tensor of order $d$ with entries
\[
[Av]_{i_1,\dots,i_d}
=
\sum_{\substack{\alpha_1,\dots,\alpha_{d-1}\\
\beta_1,\dots,\beta_{d-1}}}
C_1(i_1,\alpha_1,\beta_1)
C_2(\alpha_1,\beta_1,i_2,\alpha_2,\beta_2)
\cdots
C_d(\alpha_{d-1},\beta_{d-1},i_d),
\]
where the local contractions are defined as
\[
C_k(\alpha_k,\beta_k,i_k,\alpha_{k+1},\beta_{k+1})
=
\sum_{j_k=1}^{n_k}
G_k(\alpha_k,i_k,j_k,\alpha_{k+1})
V_k(\beta_k,j_k,\beta_{k+1}).
\]

The tensors $C_k$ can be interpreted as the TT-cores of the product $Av$. 
To be precise, the resulting representation is, in general, not in standard TT form.
Indeed, the combined indices \((\alpha_k,\beta_k)\) induce effective TT ranks $r^A_k r^v_k$
that may exceed the admissible sizes dictated by the physical dimensions \(n_k\). In such situations, a change of representation is required to restore a proper TT structure, typically by reshaping the cores and applying exact algebraic
factorizations that reorganize the indices without altering the represented tensor.

However, even if the TT ranks do not exceed the admissible sizes, explicitly
forming the new TT-cores \( C_k \) is often prohibitively expensive, as
intermediate TT ranks can grow rapidly.
Moreover, in practical applications one is typically interested in a
low-rank approximation of \( Av \), both to reduce memory requirements and to
improve computational efficiency.
In this setting, even the compression step itself may represent a significant
computational bottleneck.

This motivates the development of randomized techniques that directly compute a
compressed approximation of \( Av \).

\subsection{Khatri--Rao sketching and streaming structure}
\label{sec:khatri-rao}
\begin{figure}[ht]
    \centering
    \begin{tikzpicture}[baseline=-0.5ex, scale=0.6]
\tikzset{
  tensor/.style={circle, draw, thick, minimum size=4mm},
  green/.style={tensor, fill=green!60},
  blue/.style={tensor, fill=cyan!60},
  yellow/.style={tensor, fill=orange!80},
  wire/.style={thick},
}

\def\n{7}
\def\m{5}

\pgfmathtruncatemacro{\nhalfplusminus}{(\n+1)/2}
\pgfmathtruncatemacro{\nplus}{\n+1}
\pgfmathtruncatemacro{\nminus}{\n-1}
\pgfmathtruncatemacro{\nplusm}{\n+\m}
\pgfmathtruncatemacro{\nplusmminus}{\n+\m-1}
\pgfmathtruncatemacro{\nplusmhalf}{\n+(\m-1)/2}
\pgfmathtruncatemacro{\nplusmhalfminus}{\n-1+(\m-1)/2}
\pgfmathtruncatemacro{\nhalfminus}{(\n)/2-1}
\pgfmathtruncatemacro{\nhalfminusminus}{(\n)/2-2}
\pgfmathtruncatemacro{\nplusmhalfplus}{\n+(\m+3)/2}
\pgfmathtruncatemacro{\nhalf}{(\n-1)/2}
\pgfmathtruncatemacro{\nhalfplus}{(\n+3)/2}

\foreach \i in {1,...,\nhalf} {
  \node[green] (g\i) at (\i,0) {};
}

\foreach \i in {1,...,\nhalf} {
  \node at (\i,-1) {\tiny{$Y_{\!\scriptscriptstyle k\mkern-1.5mu+\mkern-1.5mu\i}$}};
}

\foreach \i in {1,...,\nhalfminus} {
  \node at (\n-\i,-1) {\tiny{$Y_{\!\scriptscriptstyle d\mkern-1.5mu-\mkern-1.5mu\i}$}};
}
 \node at (\n,-1) {\tiny{$Y_{\!\scriptscriptstyle d}$}};

\foreach \i in {\nhalfplus,...,\n} {
  \node[green] (g\i) at (\i,0) {};
}
\node at (\nhalf+1,0) {$\cdots$};

\foreach \i in {1,..., \nhalf} {
    \draw [wire] (g\i) -- (\i, 0.8);
}

\foreach \i in {\nhalfplus,..., \n} {
    \draw [wire] (g\i) -- (\i, 0.8);
}

\foreach \i in {1,..., \nhalf} {
    \draw [wire] (\i, -0.7) -- (\i, -0.32);
}

\foreach \i in {\nhalfplus,..., \n} {
    \draw [wire] (\i, -0.7) -- (\i, -0.32);
}

\draw [wire] (1, 0.8) -- (\n, 0.8);
\draw [wire] (\nhalfplusminus, 0.8) -- (\nhalfplusminus, 1.3);
\node[circle, draw, fill=white, inner sep=1.5pt] at ({(\n+1)/2},0.8) {};
\node at (-0.25,0) {$Y_{>k}=$};

\end{tikzpicture}
    \caption{Tensor diagram representation of $Y_{>k}$.}
    \label{fig:khatri-rao}
\end{figure}
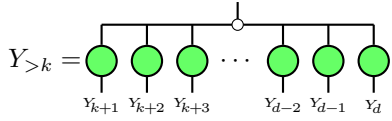

Our approach relies on structured randomized sketching based on Khatri--Rao products. 
For each mode \(i=1,\dots,d\), let \(Y_i \in \mathbb{R}^{m_i \times (s+\ell)}\) be a random sketching matrix, 
where \(s\) is the sketch size and \(\ell\) is an oversampling parameter. 
We define the nested Khatri--Rao sketches
\[
Y_{\le k} := Y_1 \odot \cdots \odot Y_k,
\qquad
Y_{>k} := Y_{k+1} \odot \cdots \odot Y_d,
\]
for \(k=1,\dots,d-1\), where \(\odot\) denotes the column-wise Kronecker (Khatri--Rao) product. In particular,
\[
[Y_{> k}]_{i_{k+1},\dots,i_d,t}
=
Y_{k+1}(i_{k+1},t)\cdot \ldots \cdot Y_d(i_d,t).
\]
This construction is illustrated in Figure~\ref{fig:khatri-rao}.
For convenience, we also define \(\bar Y_i \in \mathbb{R}^{m_i \times s}\) as the submatrix of \(Y_i\) consisting of its first \(s\) columns. 
We then define \(\bar{Y}_{\le k}\) and \(\bar{Y}_{> k}\) as the corresponding nested Khatri--Rao products.
This distinction will be useful in a subsequent construction.
Khatri--Rao sketching operators are particularly well suited to the TT format:
they exploit the separable structure of tensor trains, so that the sketching
operation decomposes along the individual modes. The alternatives are TT-Gaussian
operators, the cores of which are filled with independent Gaussian entries, and the
TT-stack constructions that interpolate between the two, trading computational
efficiency against embedding quality. Embedding guarantees for these families are
established in \cite{ahle2020oblivious,camano2025faster,cazeaux2026linear}.

Two properties of the Khatri--Rao construction matter for what follows. The first
is that the sketch factorizes into $d$ matrices $Y_k$ acting on the individual
modes, so that it can be applied one core at a time. The second is that the
column index $t$ is a \emph{batch} index rather than a summation index: in
$[Y_{>k}]_{i_{k+1},\dots,i_d,t}=Y_{k+1}(i_{k+1},t)\cdot \ldots \cdot Y_d(i_d,t)$ the same $t$
appears in every factor. Every contraction in \cref{sec:fastmatvec} is therefore
performed independently for each $t$, and the sketch dimension behaves throughout
as a batch dimension carried alongside the TT ranks. This is what allows the
matrices~$Y_k$ to be contracted into the cores of the operator without disturbing
the rest of the network.

\subsection{Streaming tensor-train approximation}
We recall the main features of the streaming tensor-train approximation algorithm \cite{kressner2023streaming}, a randomized algorithm for computing an approximation of fixed TT-rank of a tensor in tensor-train format,  as it will be central for our new method. 

The key feature of STTA is to first perform a sequence of sketching operations
associated with different matricizations of a TT tensor, and to recover the
actual approximant only after this compression phase is completed.
This approach is particularly appealing because, when forming sums or linear combinations of tensors, it suffices to sketch and compress each term
independently, combine the resulting compressed representations, and reconstruct the tensor only once at the end.
This guarantees that the algorithm can be implemented in a streamable
fashion. 

More precisely, let \(w\) be a TT tensor (here, \(w = Av\)) of which the STTA of TT-ranks \(s\) is to be computed, using the Khatri--Rao embeddings introduced above.
The STTA approximant \( \hat w \) is defined implicitly through its
\( d \)-th matricization \( \hat w_{[d]} \), which can be written as
\begin{equation*}
\begin{aligned}
\hat w_{[d]}
&=
w_{[1]} \bar{Y}_{>1}
\left( Y_{\le 1}^{\top} w_{[1]} \bar{Y}_{>1} \right)^{\dagger}
\left( Y_{\le 1}^{\top} \otimes I \right)
w_{[2]} \bar{Y}_{>2}
\left( Y_{\le 2}^{\top} w_{[2]} \bar{Y}_{>2} \right)^{\dagger}
\left( Y_{\le 2}^{\top} \otimes I \right)
\cdots \\
&\qquad
w_{[d-1]} \bar{Y}_{>d-1}
\left( Y_{\le d-1}^{\top} w_{[d-1]} \bar{Y}_{>{d-1}} \right)^{\dagger}
\left( Y_{\le d-1}^{\top} \otimes I \right)
w_{[d]}.
\end{aligned}
\end{equation*}
Here,
\[
w_{[k]} \in \mathbb{R}^{(m_1 m_2 \cdots m_k) \times (m_{k+1} \cdots m_d)}
\]
denotes the matricization of \( w \) along its first \( k \) modes.

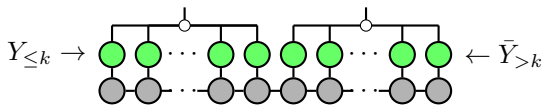
\begin{figure}[ht]
    \centering
\begin{tikzpicture}[scale=0.48]

\tikzset{
  tensor/.style={circle, draw, thick, minimum size=3mm},
  green/.style={tensor, fill=green!60},
  blue/.style={tensor, fill=cyan!60},
  yellow/.style={tensor, fill=orange!80},
  grey/.style={tensor, fill=gray!60},
  wire/.style={thick},
}

\def\n{5}
\def\m{5}
\pgfmathtruncatemacro{\nplus}{\n+1}
\pgfmathtruncatemacro{\nminus}{\n-1}
\pgfmathtruncatemacro{\nplusm}{\n+\m}
\pgfmathtruncatemacro{\nplusmminus}{\n+\m-1}
\pgfmathtruncatemacro{\nplusmhalf}{\n+(\m-1)/2}
\pgfmathtruncatemacro{\nplusmhalfminus}{\n-1+(\m-1)/2}
\pgfmathtruncatemacro{\nhalfminus}{(\n)/2-1}
\pgfmathtruncatemacro{\nplusmhalfplus}{\n+(\m+3)/2}
\pgfmathtruncatemacro{\nhalf}{(\n-1)/2}
\pgfmathtruncatemacro{\nhalfplus}{(\n+3)/2}

\foreach \i in {1,...,\nhalf} {
  \node[green] (g\i) at (\i,2) {};
}

\foreach \i in {\nhalfplus,...,\n} {
  \node[green] (g\i) at (\i,2) {};
}

\foreach \i in {\nplus,...,\nplusmhalf} {
  \node[green] (g\i) at (\i,2) {};
}
\foreach \i in {\nplusmhalfplus,...,\nplusm} {
  \node[green] (g\i) at (\i,2) {};
}

\foreach \i in {1,...,\nhalf} {
  \node[grey] (b\i) at (\i,1) {};
}

\foreach \i in {\nhalfplus,...,\n} {
  \node[grey] (b\i) at (\i,1) {};
}

\foreach \i in {\nplus,...,\nplusmhalf} {
  \node[grey] (b\i) at (\i,1) {};
}
\foreach \i in {\nplusmhalfplus,...,\nplusm} {
  \node[grey] (b\i) at (\i,1) {};
}

\foreach \i in {1,...,\nhalf} {
  \draw[wire] (g\i) -- (b\i);
}

\foreach \i in {\nhalfplus,...,\n} {
  \draw[wire] (g\i) -- (b\i);
}

\foreach \i in {\nplus,...,\nplusmhalf} {
  \draw[wire] (g\i) -- (b\i);
}

\foreach \i in {\nplusmhalfplus,...,\nplusm} {
  \draw[wire] (g\i) -- (b\i);
}

\foreach \i in {1,...,\nhalfminus} {
  \draw[wire] (b\i) -- (b\the\numexpr\i+1\relax);
}
\draw[wire] (\nhalf+0.35,1) -- (\nhalf+0.65,1);
\draw[wire] (\nhalf+1.35,1) -- (\nhalf+1.65,1);
\draw[wire] (\n+0.35,1) -- (\n+0.65,1);
\draw[wire] (\nplusmhalf+0.35,1) -- (\nplusmhalf+0.65,1);
\draw[wire] (\nplusmhalf+1.35,1) -- (\nplusmhalf+1.65,1);
\foreach \i in {\nhalfplus,...,\nminus} {
  \draw[wire] (b\i) -- (b\the\numexpr\i+1\relax);
}
\foreach \i in {\nplus,..., \nplusmhalfminus} {
  \draw[wire] (b\i) -- (b\the\numexpr\i+1\relax);
}
\foreach \i in {\nplusmhalfplus,..., \nplusmminus} {
  \draw[wire] (b\i) -- (b\the\numexpr\i+1\relax);
}



\node at (\nhalf+1,2) {$\cdots$};
\node at (\nhalf+1,1) {$\cdots$};

\node at (\nplusmhalf+1,2) {$\cdots$};
\node at (\nplusmhalf+1,1) {$\cdots$};

\draw[wire] (1,2.8) -- (\nhalf,2.8);
\foreach \i in {1,...,\nhalf} {
  \draw[wire] (g\i) -- (\i,2.8);
}
\draw[wire] (\nhalf,2.8) -- (\n,2.8);
\foreach \i in {\nhalfplus,...,\n} {
  \draw[wire] (g\i) -- (\i,2.8);
}

\draw[wire] (\nplus,2.8) -- (\nplusm,2.8);
\foreach \i in {\nplus,...,\nplusmhalf} {
  \draw[wire] (g\i) -- (\i,2.8);
}
\draw[wire] (\nhalf,2.8) -- (\n,2.8);
\foreach \i in {\nplusmhalfplus,...,\nplusm} {
  \draw[wire] (g\i) -- (\i,2.8);
}


\draw[wire] (\nhalfplus-1, 2.8) -- (\nhalfplus-1, 3.3);
\draw[wire] (\nplusmhalf+1, 2.8) -- (\nplusmhalf+1, 3.3);

\node[circle, draw, fill=white, inner sep=1.5pt] at ({(\n+1)/2},2.8) {};

\node[circle, draw, fill=white, inner sep=1.5pt] at ({(\nplus+\nplusm)/2},2.8) {};

\node at (-0.8,2) {$Y_{\leq k} \rightarrow$};
\node at (\nplusm + 1.8,2) {$\leftarrow \bar{Y}_{> k} $};

\end{tikzpicture}
\caption{Contractions required to compute $\Omega_k$.}
    \label{fig:Omega_contractions}
\end{figure}

Although the above expression may appear involved, the compression phase reduces
to the computation of two families of matrices,
specifically
\begin{equation*}
\Omega_k = Y_{\le k}^{\top} w_{[k]} \bar{Y}_{>k},
\qquad
\Psi_k = (Y_{\le k-1}^{\top} \otimes I)\, w_{[k]} \bar{Y}_{>k}.
\end{equation*}
The matrices \( \Omega_k \) are computed for \( k = 1, \ldots, d-1 \), while
\( \Psi_k \) are computed for \( k = 1, \ldots, d \).
By a slight abuse of notation, \( Y_{\le 0} \) and \( Y_{>d} \) denote identity
matrices of appropriate dimensions.
Figure~\ref{fig:Omega_contractions} shows a contraction diagram for a general $\Omega_k$. 

\begin{figure}[ht]
    \centering
    \begin{tikzpicture}[scale=0.48]

\tikzset{
  tensor/.style={circle, draw, thick, minimum size=3mm},
  green/.style={tensor, fill=green!60},
  blue/.style={tensor, fill=cyan!60},
  yellow/.style={tensor, fill=orange!80},
  grey/.style={tensor, fill=gray!60},
  wire/.style={thick},
}

\def\n{5}
\def\m{5}
\pgfmathtruncatemacro{\nplus}{\n+1}
\pgfmathtruncatemacro{\nminus}{\n-1}
\pgfmathtruncatemacro{\nplusm}{\n+\m}
\pgfmathtruncatemacro{\nplusmminus}{\n+\m-1}
\pgfmathtruncatemacro{\nplusmhalf}{\n+(\m-1)/2}
\pgfmathtruncatemacro{\nplusmhalfminus}{\n-1+(\m-1)/2}
\pgfmathtruncatemacro{\nhalfminus}{(\n)/2-1}
\pgfmathtruncatemacro{\nplusmhalfplus}{\n+(\m+3)/2}
\pgfmathtruncatemacro{\nhalf}{(\n-1)/2}
\pgfmathtruncatemacro{\nhalfplus}{(\n+3)/2}

\foreach \i in {1,...,\nhalf} {
  \node[green] (g\i) at (\i,2) {};
}

\foreach \i in {\nhalfplus,...,\nminus} {
  \node[green] (g\i) at (\i,2) {};
}

\foreach \i in {\nplus,...,\nplusmhalf} {
  \node[green] (g\i) at (\i,2) {};
}
\foreach \i in {\nplusmhalfplus,...,\nplusm} {
  \node[green] (g\i) at (\i,2) {};
}

\foreach \i in {1,...,\nhalf} {
  \node[grey] (b\i) at (\i,1) {};
}

\foreach \i in {\nhalfplus,...,\n} {
  \node[grey] (b\i) at (\i,1) {};
}

\foreach \i in {\nplus,...,\nplusmhalf} {
  \node[grey] (b\i) at (\i,1) {};
}
\foreach \i in {\nplusmhalfplus,...,\nplusm} {
  \node[grey] (b\i) at (\i,1) {};
}

\foreach \i in {1,...,\nhalf} {
  \draw[wire] (g\i) -- (b\i);
}

\foreach \i in {\nhalfplus,...,\n} {
  \draw[wire] (g\i) -- (b\i);
}

\foreach \i in {\nplus,...,\nplusmhalf} {
  \draw[wire] (g\i) -- (b\i);
}

\foreach \i in {\nplusmhalfplus,...,\nplusm} {
  \draw[wire] (g\i) -- (b\i);
}

\foreach \i in {1,...,\nhalfminus} {
  \draw[wire] (b\i) -- (b\the\numexpr\i+1\relax);
}
\draw[wire] (\nhalf+0.35,1) -- (\nhalf+0.65,1);
\draw[wire] (\nhalf+1.35,1) -- (\nhalf+1.65,1);
\draw[wire] (\n+0.35,1) -- (\n+0.65,1);
\draw[wire] (\nplusmhalf+0.35,1) -- (\nplusmhalf+0.65,1);
\draw[wire] (\nplusmhalf+1.35,1) -- (\nplusmhalf+1.65,1);
\foreach \i in {\nhalfplus,...,\nminus} {
  \draw[wire] (b\i) -- (b\the\numexpr\i+1\relax);
}
\foreach \i in {\nplus,..., \nplusmhalfminus} {
  \draw[wire] (b\i) -- (b\the\numexpr\i+1\relax);
}
\foreach \i in {\nplusmhalfplus,..., \nplusmminus} {
  \draw[wire] (b\i) -- (b\the\numexpr\i+1\relax);
}

\node at (\nhalf+1,2) {$\cdots$};
\node at (\nhalf+1,1) {$\cdots$};
\node at (\nplusmhalf+1,2) {$\cdots$};
\node at (\nplusmhalf+1,1) {$\cdots$};



\draw[wire] (1,2.8) -- (\nhalf,2.8);
\foreach \i in {1,...,\nhalf} {
  \draw[wire] (g\i) -- (\i,2.8);
}

\foreach \i in {\nhalfplus,...,\nminus} {
  \draw[wire] (g\i) -- (\i,2.8);
}

\draw[wire] (\nplus,2.8) -- (\nplusm,2.8);
\foreach \i in {\nplus,...,\nplusmhalf} {
  \draw[wire] (g\i) -- (\i,2.8);
}
\draw[wire] (\nhalf,2.8) -- (\nminus,2.8);
\foreach \i in {\nplusmhalfplus,...,\nplusm} {
  \draw[wire] (g\i) -- (\i,2.8);
}

\draw[wire] (\nhalfplus-1, 2.8) -- (\nhalfplus-1, 3.3);
\draw[wire] (\nplusmhalf+1, 2.8) -- (\nplusmhalf+1, 3.3);

\node[circle, draw, fill=white, inner sep=1.5pt] at ({(\n+1)/2},2.8) {};

\node[circle, draw, fill=white, inner sep=1.5pt] at ({(\nplus+\nplusm)/2},2.8) {};

\node at (-0.8,2) {$Y_{\leq k-1} \rightarrow$};
\node at (\nplusm + 1.8,2) {$\leftarrow \bar{Y}_{> k} $};

\end{tikzpicture}
    \caption{Contractions required to compute $\Psi_k$.}
    \label{fig:Psi_contractions}
\end{figure}
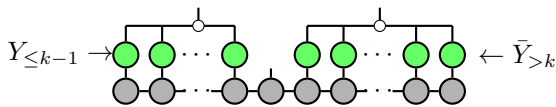

Similarly, Figure~\ref{fig:Psi_contractions} represents the contraction diagram for $\Psi_k$.

To recover an actual tensor-train approximation, we still need to compute the TT-cores. 
This operation is not linear in the tensor, as it involves the pseudo-inverses. Consequently, if one wishes to compute a sum of tensors, say \(w^{(1)}, \dots, w^{(q)}\),
it suffices to sum the corresponding \(\Omega_k^{(p)}\) and \(\Psi_k^{(p)}\), for $p=1,\dots, q$ first,
and then form the TT approximation from these combined matrices.
The TT-cores of the approximation (or, more precisely, the corresponding matricizations) are then given by
\begin{equation*}
    W_k =
    \begin{cases}
        \Psi_1, & k = 1,\\[2mm]
        \Omega_{k-1}^\dagger \Psi_k, & k > 1.
    \end{cases}
\end{equation*}
From an implementation perspective, the pseudoinverse can be highly ill-conditioned, potentially leading to numerical issues. 
For this reason, following~\cite{nakatsukasa2020fast}, 
we apply an $\epsilon$-truncation of the pseudoinverse in our experiments, setting 
$ \epsilon = u_{\text{mach}} \, \|\Omega_{k-1}\|_2$, where $u_{\text{mach}}$ denotes the machine precision.  This solely requires an SVD of the small $(s+\ell)\times s$ matrix $\Omega_{k-1}$. 
We then discard all singular values of $\Omega_{k-1}$ smaller than $\epsilon$.

Another useful tool in the tensor-train framework is that of an
$( \varepsilon, \delta, m)$-oblivious subspace embedding (OSE). Formally, a linear operator $S \in \mathbb{R}^{s \times M}$ is an
$(\varepsilon, \delta, m)$-OSE if, for any $m$-dimensional subspace
$\mathcal{X} \subseteq \mathbb{R}^M$, with probability at least $1-\delta$ it holds that for all $x \in \mathcal{X}$,
\begin{equation*}
    (1-\varepsilon)\|x\|_2^2 \le \|S x\|_2^2 \le (1+\varepsilon)\|x\|_2^2.
\end{equation*}
When working with tensors, natural choices for such embeddings include Khatri--Rao operators, TT-structured maps (e.g., TT-Gaussian maps), and their corresponding block variants. In this work, we focus on Khatri--Rao sketching; however, the analysis extends directly to the other structured embeddings mentioned above.

\section{Fast TT matrix--vector multiplication}
\label{sec:fastmatvec}

We now describe how the compressed matrices $\Omega_k$ and $\Psi_k$ of
\cref{sec:preliminaries}, from which the TT approximation of $w=Av$ is recovered,
can be assembled without ever forming the whole product. The key observation is that, at each site, the core of the operator, the core of the vector, and the sketching matrix can be contracted in different orders. The choice of order determines which intermediate tensors need to be formed.

\subsection{Three contraction orderings}
\label{sec:orders}

Write $G_k\in\mathbb{R}^{r_{k-1}^A\times m_k\times n_k\times r_k^A}$ for the
$k$-th core of $A$, $V_k\in\mathbb{R}^{r_{k-1}^v\times n_k\times r_k^v}$ for the
$k$-th core of $v$, and $Y_k\in\mathbb{R}^{m_k\times(s+\ell)}$ for the $k$-th
sketching matrix. Every quantity we need is obtained by contracting these
tensors, and three contraction orderings are of interest.

The \emph{sequential ordering} first contracts $G_k$ with $V_k$. The result is the
$k$-th core of~$Av$, of size
$r_{k-1}^A r_{k-1}^v\times m_k\times r_k^A r_k^v$.
Forming this core costs
$\mathcal{O}\!\left(m_k n_k r_{k-1}^A r_{k-1}^vr_k^A r_k^v\right)$,
which is quadratic in the TT ranks of both $A$ and $v$. Once the core has been
formed, the compression proceeds exactly as in STTA applied to a given TT
tensor. Only one core of $Av$ needs to be stored at a time, so the full product
need not be formed explicitly, but each of its cores must still be constructed.

The \emph{fused ordering} avoids explicitly forming the cores of $Av$. Instead,
$V_k$ is contracted with the partial contraction arriving from the previous
site, after which~$G_k$ is applied to the result and $Y_k$ is applied last.

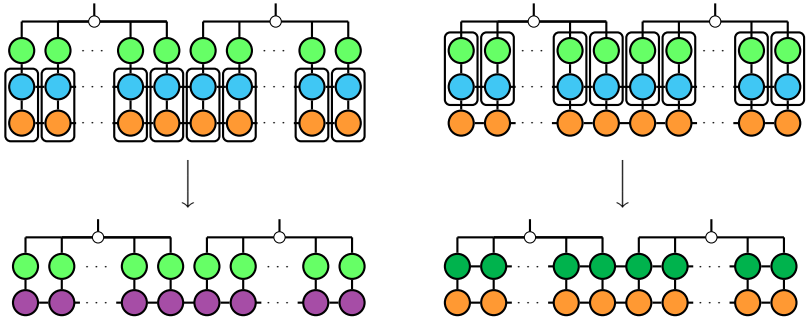
\begin{figure}[ht]
    \centering
\begin{tikzpicture}[baseline=-0.5ex, scale=0.48]

\tikzset{
  tensor/.style={circle, draw, thick, minimum size=3mm},
  green/.style={tensor, fill=green!60},
  blue/.style={tensor, fill=cyan!60},
  yellow/.style={tensor, fill=orange!80},
  wire/.style={thick},
}

\def\n{5}
\def\m{5}
\pgfmathtruncatemacro{\nplus}{\n+1}
\pgfmathtruncatemacro{\nminus}{\n-1}
\pgfmathtruncatemacro{\nplusm}{\n+\m}
\pgfmathtruncatemacro{\nplusmminus}{\n+\m-1}
\pgfmathtruncatemacro{\nplusmhalf}{\n+(\m-1)/2}
\pgfmathtruncatemacro{\nplusmhalfminus}{\n-1+(\m-1)/2}
\pgfmathtruncatemacro{\nhalfminus}{(\n)/2-1}
\pgfmathtruncatemacro{\nplusmhalfplus}{\n+(\m+3)/2}
\pgfmathtruncatemacro{\nhalf}{(\n-1)/2}
\pgfmathtruncatemacro{\nhalfplus}{(\n+3)/2}

\foreach \i in {1,...,\nhalf} {
  \node[green] (g\i) at (\i,2) {};
}

\foreach \i in {\nhalfplus,...,\n} {
  \node[green] (g\i) at (\i,2) {};
}

\foreach \i in {\nplus,...,\nplusmhalf} {
  \node[green] (g\i) at (\i,2) {};
}
\foreach \i in {\nplusmhalfplus,...,\nplusm} {
  \node[green] (g\i) at (\i,2) {};
}

\foreach \i in {1,...,\nhalf} {
  \node[blue] (b\i) at (\i,1) {};
}

\foreach \i in {\nhalfplus,...,\n} {
  \node[blue] (b\i) at (\i,1) {};
}

\foreach \i in {\nplus,...,\nplusmhalf} {
  \node[blue] (b\i) at (\i,1) {};
}
\foreach \i in {\nplusmhalfplus,...,\nplusm} {
  \node[blue] (b\i) at (\i,1) {};
}

\foreach \i in {1,...,\nhalf} {
  \node[yellow] (y\i) at (\i,0) {};
}

\foreach \i in {\nhalfplus,...,\n} {
  \node[yellow] (y\i) at (\i,0) {};
}

\foreach \i in {\nplus,...,\nplusmhalf} {
  \node[yellow] (y\i) at (\i,0) {};
}
\foreach \i in {\nplusmhalfplus,...,\nplusm} {
  \node[yellow] (y\i) at (\i,0) {};
}

\foreach \i in {1,...,\nhalf} {
  \draw[wire] (g\i) -- (b\i);
  \draw[wire] (b\i) -- (y\i);
}

\foreach \i in {\nhalfplus,...,\n} {
  \draw[wire] (g\i) -- (b\i);
  \draw[wire] (b\i) -- (y\i);
}

\foreach \i in {\nplus,...,\nplusmhalf} {
  \draw[wire] (g\i) -- (b\i);
  \draw[wire] (b\i) -- (y\i);
}

\foreach \i in {\nplusmhalfplus,...,\nplusm} {
  \draw[wire] (g\i) -- (b\i);
  \draw[wire] (b\i) -- (y\i);
}

\foreach \i in {1,...,\nhalfminus} {
  \draw[wire] (b\i) -- (b\the\numexpr\i+1\relax);
}
\draw[wire] (\nhalf+0.35,1) -- (\nhalf+0.5,1);
\draw[wire] (\nhalf+1.5,1) -- (\nhalf+1.65,1);
\draw[wire] (\n+0.35,1) -- (\n+0.65,1);
\draw[wire] (\nplusmhalf+0.35,1) -- (\nplusmhalf+0.5,1);
\draw[wire] (\nplusmhalf+1.5,1) -- (\nplusmhalf+1.65,1);
\foreach \i in {\nhalfplus,...,\nminus} {
  \draw[wire] (b\i) -- (b\the\numexpr\i+1\relax);
}
\foreach \i in {\nplus,..., \nplusmhalfminus} {
  \draw[wire] (b\i) -- (b\the\numexpr\i+1\relax);
}
\foreach \i in {\nplusmhalfplus,..., \nplusmminus} {
  \draw[wire] (b\i) -- (b\the\numexpr\i+1\relax);
}

\foreach \i in {1,...,\nhalfminus} {
  \draw[wire] (y\i) -- (y\the\numexpr\i+1\relax);
}
\draw[wire] (\nhalf+0.35,0) -- (\nhalf+0.5,0);
\draw[wire] (\nhalf+1.5,0) -- (\nhalf+1.65,0);
\draw[wire] (\n+0.35,0) -- (\n+0.65,0);
\foreach \i in {\nplus,..., \nplusmhalfminus} {
  \draw[wire] (y\i) -- (y\the\numexpr\i+1\relax);
}
\draw[wire] (\nplusmhalf+0.35,0) -- (\nplusmhalf+0.5,0);
\draw[wire] (\nplusmhalf+1.5,0) -- (\nplusmhalf+1.65,0);
\foreach \i in {\nhalfplus,...,\nminus} {
  \draw[wire] (y\i) -- (y\the\numexpr\i+1\relax);
}
\foreach \i in {\nplusmhalfplus,..., \nplusmminus} {
  \draw[wire] (y\i) -- (y\the\numexpr\i+1\relax);
}

\node at (\nhalf+1,2) {\tiny $\cdots$};
\node at (\nhalf+1,1) {\tiny $\cdots$};
\node at (\nhalf+1,0) {\tiny$\cdots$};

\node at (\nplusmhalf+1,2) {\tiny$\cdots$};
\node at (\nplusmhalf+1,1) {\tiny $\cdots$};
\node at (\nplusmhalf+1,0) {\tiny $\cdots$};

\draw[wire] (1,2.8) -- (\nhalf,2.8);
\foreach \i in {1,...,\nhalf} {
  \draw[wire] (g\i) -- (\i,2.8);
}
\draw[wire] (\nhalf,2.8) -- (\n,2.8);
\foreach \i in {\nhalfplus,...,\n} {
  \draw[wire] (g\i) -- (\i,2.8);
}

\draw[wire] (\nplus,2.8) -- (\nplusm,2.8);
\foreach \i in {\nplus,...,\nplusmhalf} {
  \draw[wire] (g\i) -- (\i,2.8);
}
\draw[wire] (\nhalf,2.8) -- (\n,2.8);
\foreach \i in {\nplusmhalfplus,...,\nplusm} {
  \draw[wire] (g\i) -- (\i,2.8);
}

\draw[wire] (\nhalfplus-1, 2.8) -- (\nhalfplus-1, 3.3);
\draw[wire] (\nplusmhalf+1, 2.8) -- (\nplusmhalf+1, 3.3);
\node[circle, draw, fill=white, inner sep=1.5pt] at ({(\n+1)/2},2.8) {};

\node[circle, draw, fill=white, inner sep=1.5pt] at ({(\nplus+\nplusm)/2},2.8) {};

\draw[draw=black, fill=none, rounded corners=2pt, thick] (1.45,1.5) rectangle (0.55,-0.5);
\draw[draw=black, fill=none, rounded corners=2pt, thick] (2.45,1.5) rectangle (1.55,-0.5);
\draw[draw=black, fill=none, rounded corners=2pt, thick] (4.45,1.5) rectangle (3.55,-0.5);
\draw[draw=black, fill=none, rounded corners=2pt, thick] (5.45,1.5) rectangle (4.55,-0.5);
\draw[draw=black, fill=none, rounded corners=2pt, thick] (6.45,1.5) rectangle (5.55,-0.5);
\draw[draw=black, fill=none, rounded corners=2pt, thick] (7.45,1.5) rectangle (6.55,-0.5);
\draw[draw=black, fill=none, rounded corners=2pt, thick] (9.45,1.5) rectangle (8.55,-0.5);
\draw[draw=black, fill=none, rounded corners=2pt, thick] (10.45,1.5) rectangle (9.55,-0.5);
\end{tikzpicture}
\hspace{8mm}
\begin{tikzpicture}[baseline=-0.5ex, scale=0.48]
\tikzset{
  tensor/.style={circle, draw, thick, minimum size=3mm},
  green/.style={tensor, fill=green!60},
  blue/.style={tensor, fill=cyan!60},
  yellow/.style={tensor, fill=orange!80},
  wire/.style={thick},
}

\def\n{5}
\def\m{5}
\pgfmathtruncatemacro{\nplus}{\n+1}
\pgfmathtruncatemacro{\nminus}{\n-1}
\pgfmathtruncatemacro{\nplusm}{\n+\m}
\pgfmathtruncatemacro{\nplusmminus}{\n+\m-1}
\pgfmathtruncatemacro{\nplusmhalf}{\n+(\m-1)/2}
\pgfmathtruncatemacro{\nplusmhalfminus}{\n-1+(\m-1)/2}
\pgfmathtruncatemacro{\nhalfminus}{(\n)/2-1}
\pgfmathtruncatemacro{\nplusmhalfplus}{\n+(\m+3)/2}
\pgfmathtruncatemacro{\nhalf}{(\n-1)/2}
\pgfmathtruncatemacro{\nhalfplus}{(\n+3)/2}

\foreach \i in {1,...,\nhalf} {
  \node[green] (g\i) at (\i,2) {};
}

\foreach \i in {\nhalfplus,...,\n} {
  \node[green] (g\i) at (\i,2) {};
}

\foreach \i in {\nplus,...,\nplusmhalf} {
  \node[green] (g\i) at (\i,2) {};
}
\foreach \i in {\nplusmhalfplus,...,\nplusm} {
  \node[green] (g\i) at (\i,2) {};
}

\foreach \i in {1,...,\nhalf} {
  \node[blue] (b\i) at (\i,1) {};
}

\foreach \i in {\nhalfplus,...,\n} {
  \node[blue] (b\i) at (\i,1) {};
}

\foreach \i in {\nplus,...,\nplusmhalf} {
  \node[blue] (b\i) at (\i,1) {};
}
\foreach \i in {\nplusmhalfplus,...,\nplusm} {
  \node[blue] (b\i) at (\i,1) {};
}

\foreach \i in {1,...,\nhalf} {
  \node[yellow] (y\i) at (\i,0) {};
}

\foreach \i in {\nhalfplus,...,\n} {
  \node[yellow] (y\i) at (\i,0) {};
}

\foreach \i in {\nplus,...,\nplusmhalf} {
  \node[yellow] (y\i) at (\i,0) {};
}
\foreach \i in {\nplusmhalfplus,...,\nplusm} {
  \node[yellow] (y\i) at (\i,0) {};
}

\foreach \i in {1,...,\nhalf} {
  \draw[wire] (g\i) -- (b\i);
  \draw[wire] (b\i) -- (y\i);
}

\foreach \i in {\nhalfplus,...,\n} {
  \draw[wire] (g\i) -- (b\i);
  \draw[wire] (b\i) -- (y\i);
}

\foreach \i in {\nplus,...,\nplusmhalf} {
  \draw[wire] (g\i) -- (b\i);
  \draw[wire] (b\i) -- (y\i);
}

\foreach \i in {\nplusmhalfplus,...,\nplusm} {
  \draw[wire] (g\i) -- (b\i);
  \draw[wire] (b\i) -- (y\i);
}

\foreach \i in {1,...,\nhalfminus} {
  \draw[wire] (b\i) -- (b\the\numexpr\i+1\relax);
}
\draw[wire] (\nhalf+0.35,1) -- (\nhalf+0.5,1);
\draw[wire] (\nhalf+1.5,1) -- (\nhalf+1.65,1);
\draw[wire] (\n+0.35,1) -- (\n+0.65,1);
\draw[wire] (\nplusmhalf+0.35,1) -- (\nplusmhalf+0.5,1);
\draw[wire] (\nplusmhalf+1.5,1) -- (\nplusmhalf+1.65,1);
\foreach \i in {\nhalfplus,...,\nminus} {
  \draw[wire] (b\i) -- (b\the\numexpr\i+1\relax);
}
\foreach \i in {\nplus,..., \nplusmhalfminus} {
  \draw[wire] (b\i) -- (b\the\numexpr\i+1\relax);
}
\foreach \i in {\nplusmhalfplus,..., \nplusmminus} {
  \draw[wire] (b\i) -- (b\the\numexpr\i+1\relax);
}

\foreach \i in {1,...,\nhalfminus} {
  \draw[wire] (y\i) -- (y\the\numexpr\i+1\relax);
}
\draw[wire] (\nhalf+0.35,0) -- (\nhalf+0.5,0);
\draw[wire] (\nhalf+1.5,0) -- (\nhalf+1.65,0);
\draw[wire] (\n+0.35,0) -- (\n+0.65,0);
\foreach \i in {\nplus,..., \nplusmhalfminus} {
  \draw[wire] (y\i) -- (y\the\numexpr\i+1\relax);
}
\draw[wire] (\nplusmhalf+0.35,0) -- (\nplusmhalf+0.5,0);
\draw[wire] (\nplusmhalf+1.5,0) -- (\nplusmhalf+1.65,0);
\foreach \i in {\nhalfplus,...,\nminus} {
  \draw[wire] (y\i) -- (y\the\numexpr\i+1\relax);
}
\foreach \i in {\nplusmhalfplus,..., \nplusmminus} {
  \draw[wire] (y\i) -- (y\the\numexpr\i+1\relax);
}

\node at (\nhalf+1,2) {\tiny $\cdots$};
\node at (\nhalf+1,1) {\tiny $\cdots$};
\node at (\nhalf+1,0) {\tiny$\cdots$};

\node at (\nplusmhalf+1,2) {\tiny$\cdots$};
\node at (\nplusmhalf+1,1) {\tiny $\cdots$};
\node at (\nplusmhalf+1,0) {\tiny $\cdots$};



\draw[wire] (1,2.8) -- (\nhalf,2.8);
\foreach \i in {1,...,\nhalf} {
  \draw[wire] (g\i) -- (\i,2.8);
}
\draw[wire] (\nhalf,2.8) -- (\n,2.8);
\foreach \i in {\nhalfplus,...,\n} {
  \draw[wire] (g\i) -- (\i,2.8);
}

\draw[wire] (\nplus,2.8) -- (\nplusm,2.8);
\foreach \i in {\nplus,...,\nplusmhalf} {
  \draw[wire] (g\i) -- (\i,2.8);
}
\draw[wire] (\nhalf,2.8) -- (\n,2.8);
\foreach \i in {\nplusmhalfplus,...,\nplusm} {
  \draw[wire] (g\i) -- (\i,2.8);
}

\draw[wire] (\nhalfplus-1, 2.8) -- (\nhalfplus-1, 3.3);
\draw[wire] (\nplusmhalf+1, 2.8) -- (\nplusmhalf+1, 3.3);
\node[circle, draw, fill=white, inner sep=1.5pt] at ({(\n+1)/2},2.8) {};

\node[circle, draw, fill=white, inner sep=1.5pt] at ({(\nplus+\nplusm)/2},2.8) {};

\draw[draw=black, fill=none, rounded corners=2pt, thick] (1.45,2.5) rectangle (0.55,0.5);
\draw[draw=black, fill=none, rounded corners=2pt, thick] (2.45,2.5) rectangle (1.55,0.5);
\draw[draw=black, fill=none, rounded corners=2pt, thick] (4.45,2.5) rectangle (3.55,0.5);
\draw[draw=black, fill=none, rounded corners=2pt, thick] (5.45,2.5) rectangle (4.55,0.5);
\draw[draw=black, fill=none, rounded corners=2pt, thick] (6.45,2.5) rectangle (5.55,0.5);
\draw[draw=black, fill=none, rounded corners=2pt, thick] (7.45,2.5) rectangle (6.55,0.5);
\draw[draw=black, fill=none, rounded corners=2pt, thick] (9.45,2.5) rectangle (8.55,0.5);
\draw[draw=black, fill=none, rounded corners=2pt, thick] (10.45,2.5) rectangle (9.55,0.5);

\end{tikzpicture}
\vspace{2mm}\\
\hspace{-1mm}$\Big\downarrow$
\hspace{54mm}$\Big\downarrow$\vspace{1mm}\\
    \begin{tikzpicture}[baseline=-0.5ex, scale=0.48]
\tikzset{
  tensor/.style={circle, draw, thick, minimum size=3mm},
  green/.style={tensor, fill=green!60},
  blue/.style={tensor, fill=cyan!60},
  yellow/.style={tensor, fill=violet!70!white},
  avg/.style={tensor, fill=green!60},
  wire/.style={thick},
}

\def\n{5}
\def\m{5}
\pgfmathtruncatemacro{\nplus}{\n+1}
\pgfmathtruncatemacro{\nminus}{\n-1}
\pgfmathtruncatemacro{\nplusm}{\n+\m}
\pgfmathtruncatemacro{\nplusmminus}{\n+\m-1}
\pgfmathtruncatemacro{\nplusmhalf}{\n+(\m-1)/2}
\pgfmathtruncatemacro{\nplusmhalfminus}{\n-1+(\m-1)/2}
\pgfmathtruncatemacro{\nhalfminus}{(\n)/2-1}
\pgfmathtruncatemacro{\nplusmhalfplus}{\n+(\m+3)/2}
\pgfmathtruncatemacro{\nhalf}{(\n-1)/2}
\pgfmathtruncatemacro{\nhalfplus}{(\n+3)/2}

\foreach \i in {1,...,\nhalf} {
  \node[avg] (b\i) at (\i,1) {};
}

\foreach \i in {\nhalfplus,...,\n} {
  \node[avg] (b\i) at (\i,1) {};
}

\foreach \i in {\nplus,...,\nplusmhalf} {
  \node[avg] (b\i) at (\i,1) {};
}
\foreach \i in {\nplusmhalfplus,...,\nplusm} {
  \node[avg] (b\i) at (\i,1) {};
}

\foreach \i in {1,...,\nhalf} {
  \node[yellow] (y\i) at (\i,0) {};
}

\foreach \i in {\nhalfplus,...,\n} {
  \node[yellow] (y\i) at (\i,0) {};
}

\foreach \i in {\nplus,...,\nplusmhalf} {
  \node[yellow] (y\i) at (\i,0) {};
}
\foreach \i in {\nplusmhalfplus,...,\nplusm} {
  \node[yellow] (y\i) at (\i,0) {};
}

\foreach \i in {1,...,\nhalf} {
  \draw[wire] (b\i) -- (y\i);
}

\foreach \i in {\nhalfplus,...,\n} {
  \draw[wire] (b\i) -- (y\i);
}

\foreach \i in {\nplus,...,\nplusmhalf} {
  \draw[wire] (b\i) -- (y\i);
}

\foreach \i in {\nplusmhalfplus,...,\nplusm} {
  \draw[wire] (b\i) -- (y\i);
}

\foreach \i in {1,...,\nhalfminus} {
  \draw[wire] (y\i) -- (y\the\numexpr\i+1\relax);
}
\draw[wire] (\nhalf+0.35,0) -- (\nhalf+0.5,0);
\draw[wire] (\nhalf+1.5,0) -- (\nhalf+1.65,0);
\draw[wire] (\n+0.35,0) -- (\n+0.65,0);
\foreach \i in {\nplus,..., \nplusmhalfminus} {
  \draw[wire] (y\i) -- (y\the\numexpr\i+1\relax);
}
\draw[wire] (\nplusmhalf+0.35,0) -- (\nplusmhalf+0.5,0);
\draw[wire] (\nplusmhalf+1.5,0) -- (\nplusmhalf+1.65,0);
\foreach \i in {\nhalfplus,...,\nminus} {
  \draw[wire] (y\i) -- (y\the\numexpr\i+1\relax);
}
\foreach \i in {\nplusmhalfplus,..., \nplusmminus} {
  \draw[wire] (y\i) -- (y\the\numexpr\i+1\relax);
}

\node at (\nhalf+1,1) {\tiny $\cdots$};
\node at (\nhalf+1,0) {\tiny$\cdots$};

\node at (\nplusmhalf+1,1) {\tiny $\cdots$};
\node at (\nplusmhalf+1,0) {\tiny $\cdots$};



\draw[wire] (1,1.8) -- (\nhalf,1.8);
\foreach \i in {1,...,\nhalf} {
  \draw[wire] (b\i) -- (\i,1.8);
}
\draw[wire] (\nhalf,1.8) -- (\n,1.8);
\foreach \i in {\nhalfplus,...,\n} {
  \draw[wire] (b\i) -- (\i,1.8);
}

\draw[wire] (\nplus,1.8) -- (\nplusm,1.8);
\foreach \i in {\nplus,...,\nplusmhalf} {
  \draw[wire] (b\i) -- (\i,1.8);
}
\draw[wire] (\nhalf,1.8) -- (\n,1.8);
\foreach \i in {\nplusmhalfplus,...,\nplusm} {
  \draw[wire] (b\i) -- (\i,1.8);
}

\draw[wire] (\nhalfplus-1, 1.8) -- (\nhalfplus-1, 2.3);
\draw[wire] (\nplusmhalf+1, 1.8) -- (\nplusmhalf+1, 2.3);
\node[circle, draw, fill=white, inner sep=1.5pt] at ({(\n+1)/2},1.8) {};

\node[circle, draw, fill=white, inner sep=1.5pt] at ({(\nplus+\nplusm)/2},1.8) {};

\end{tikzpicture}
\hspace{8mm}
    \begin{tikzpicture}[baseline=-0.5ex, scale=0.48]
\tikzset{
  tensor/.style={circle, draw, thick, minimum size=3mm},
  green/.style={tensor, fill=green!60},
  blue/.style={tensor, fill=cyan!60},
  yellow/.style={tensor, fill=orange!80},
  avg/.style={tensor, fill=green!70!blue},
  wire/.style={thick},
}

\def\n{5}
\def\m{5}
\pgfmathtruncatemacro{\nplus}{\n+1}
\pgfmathtruncatemacro{\nminus}{\n-1}
\pgfmathtruncatemacro{\nplusm}{\n+\m}
\pgfmathtruncatemacro{\nplusmminus}{\n+\m-1}
\pgfmathtruncatemacro{\nplusmhalf}{\n+(\m-1)/2}
\pgfmathtruncatemacro{\nplusmhalfminus}{\n-1+(\m-1)/2}
\pgfmathtruncatemacro{\nhalfminus}{(\n)/2-1}
\pgfmathtruncatemacro{\nplusmhalfplus}{\n+(\m+3)/2}
\pgfmathtruncatemacro{\nhalf}{(\n-1)/2}
\pgfmathtruncatemacro{\nhalfplus}{(\n+3)/2}

\foreach \i in {1,...,\nhalf} {
  \node[avg] (b\i) at (\i,1) {};
}

\foreach \i in {\nhalfplus,...,\n} {
  \node[avg] (b\i) at (\i,1) {};
}

\foreach \i in {\nplus,...,\nplusmhalf} {
  \node[avg] (b\i) at (\i,1) {};
}
\foreach \i in {\nplusmhalfplus,...,\nplusm} {
  \node[avg] (b\i) at (\i,1) {};
}

\foreach \i in {1,...,\nhalf} {
  \node[yellow] (y\i) at (\i,0) {};
}

\foreach \i in {\nhalfplus,...,\n} {
  \node[yellow] (y\i) at (\i,0) {};
}

\foreach \i in {\nplus,...,\nplusmhalf} {
  \node[yellow] (y\i) at (\i,0) {};
}
\foreach \i in {\nplusmhalfplus,...,\nplusm} {
  \node[yellow] (y\i) at (\i,0) {};
}

\foreach \i in {1,...,\nhalf} {
  \draw[wire] (b\i) -- (y\i);
}

\foreach \i in {\nhalfplus,...,\n} {
  \draw[wire] (b\i) -- (y\i);
}

\foreach \i in {\nplus,...,\nplusmhalf} {
  \draw[wire] (b\i) -- (y\i);
}

\foreach \i in {\nplusmhalfplus,...,\nplusm} {
  \draw[wire] (b\i) -- (y\i);
}

\foreach \i in {1,...,\nhalfminus} {
  \draw[wire] (b\i) -- (b\the\numexpr\i+1\relax);
}
\draw[wire] (\nhalf+0.35,1) -- (\nhalf+0.5,1);
\draw[wire] (\nhalf+1.5,1) -- (\nhalf+1.65,1);
\draw[wire] (\n+0.35,1) -- (\n+0.65,1);
\draw[wire] (\nplusmhalf+0.35,1) -- (\nplusmhalf+0.5,1);
\draw[wire] (\nplusmhalf+1.5,1) -- (\nplusmhalf+1.65,1);
\foreach \i in {\nhalfplus,...,\nminus} {
  \draw[wire] (b\i) -- (b\the\numexpr\i+1\relax);
}
\foreach \i in {\nplus,..., \nplusmhalfminus} {
  \draw[wire] (b\i) -- (b\the\numexpr\i+1\relax);
}
\foreach \i in {\nplusmhalfplus,..., \nplusmminus} {
  \draw[wire] (b\i) -- (b\the\numexpr\i+1\relax);
}

\foreach \i in {1,...,\nhalfminus} {
  \draw[wire] (y\i) -- (y\the\numexpr\i+1\relax);
}
\draw[wire] (\nhalf+0.35,0) -- (\nhalf+0.5,0);
\draw[wire] (\nhalf+1.5,0) -- (\nhalf+1.65,0);
\draw[wire] (\n+0.35,0) -- (\n+0.65,0);
\foreach \i in {\nplus,..., \nplusmhalfminus} {
  \draw[wire] (y\i) -- (y\the\numexpr\i+1\relax);
}
\draw[wire] (\nplusmhalf+0.35,0) -- (\nplusmhalf+0.5,0);
\draw[wire] (\nplusmhalf+1.5,0) -- (\nplusmhalf+1.65,0);
\foreach \i in {\nhalfplus,...,\nminus} {
  \draw[wire] (y\i) -- (y\the\numexpr\i+1\relax);
}
\foreach \i in {\nplusmhalfplus,..., \nplusmminus} {
  \draw[wire] (y\i) -- (y\the\numexpr\i+1\relax);
}

\node at (\nhalf+1,1) {\tiny $\cdots$};
\node at (\nhalf+1,0) {\tiny$\cdots$};

\node at (\nplusmhalf+1,1) {\tiny $\cdots$};
\node at (\nplusmhalf+1,0) {\tiny $\cdots$};



\draw[wire] (1,1.8) -- (\nhalf,1.8);
\foreach \i in {1,...,\nhalf} {
  \draw[wire] (b\i) -- (\i,1.8);
}
\draw[wire] (\nhalf,1.8) -- (\n,1.8);
\foreach \i in {\nhalfplus,...,\n} {
  \draw[wire] (b\i) -- (\i,1.8);
}

\draw[wire] (\nplus,1.8) -- (\nplusm,1.8);
\foreach \i in {\nplus,...,\nplusmhalf} {
  \draw[wire] (b\i) -- (\i,1.8);
}
\draw[wire] (\nhalf,1.8) -- (\n,1.8);
\foreach \i in {\nplusmhalfplus,...,\nplusm} {
  \draw[wire] (b\i) -- (\i,1.8);
}

\draw[wire] (\nhalfplus-1, 1.8) -- (\nhalfplus-1, 2.3);
\draw[wire] (\nplusmhalf+1, 1.8) -- (\nplusmhalf+1, 2.3);
\node[circle, draw, fill=white, inner sep=1.5pt] at ({(\n+1)/2},1.8) {};

\node[circle, draw, fill=white, inner sep=1.5pt] at ({(\nplus+\nplusm)/2},1.8) {};

\end{tikzpicture}
\caption{Diagrammatic representation of two contraction strategies for computing the compressed matrices $\Omega_k$. The rectangular boxes highlight the contractions performed first. Left: The cores of the operator and the vector are contracted first to form the cores of $Av$, which are then sketched, corresponding to first forming $Av$ and subsequently applying standard STTA compression. Right: The pre-sketched approach, where the sketch is absorbed into the operator cores beforehand, leaving a much smaller object to be contracted with the vector.}
    \label{fig:matvec_Omega_contractions}
\end{figure}

The \emph{pre-sketched ordering} goes one step further by contracting each $Y_k$
with the corresponding $G_k$ before the remaining contractions are performed.
Specifically, we form
$$
  \widetilde G_k(\alpha,j,\alpha',t)
  \;=\sum_{i=1}^{m_k}Y_k(i,t)\,G_k(\alpha,i,j,\alpha'),
$$
once and for all. The resulting tensor has size
$r_{k-1}^A\times n_k\times r_k^A\times(s+\ell)$, with the row mode $m_k$
replaced by the sketch dimension $s+\ell$. The remaining contractions are then
the same as in the fused order, but with $\widetilde G_k$ in place of $G_k$;
consequently, the row modes of the operator do not enter these contractions.
The contraction defining~$\widetilde G_k$ is illustrated in
Figure~\ref{fig:matvec_Omega_contractions}.

The three orders contract the same network and therefore return the same
$\Omega_k$ and~$\Psi_k$, hence the same approximant, namely the STTA
approximation of $Av$ obtained with the sketching matrices $Y_k$. In particular, the
quasi-optimality guarantees established for STTA apply verbatim to all three and the three differ only in cost. They also inherit the linearity of $\Omega_k$ and $\Psi_k$ in the tensor, allowing the terms in a linear combination $\sum_p c_p Av^{(p)}$ to be compressed separately and combined to obtain a low-rank approximation of the sum, as used in \cref{sec:exp-stream}.

\subsection{Organizing the sweeps}
\label{sec:sweeps}

Whichever ordering is chosen, the contractions are arranged so that intermediate
results are shared. They start at the peripheral cores and proceed inwards: the
left contractions from $k=1$ towards $k=d$, the right contractions from $k=d$
backwards. The two sequences are independent and may be computed in parallel.

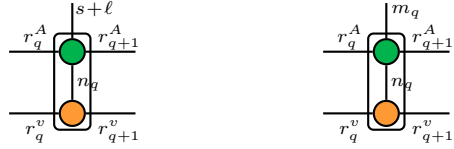
\begin{figure}[ht]
\centering
\begin{tikzpicture}[baseline=-0.5ex, scale=0.48]
\tikzset{
  tensor/.style={circle, draw, thick, minimum size=3mm},
  green/.style={tensor, fill=green!60},
  blue/.style={tensor, fill=cyan!60},
  yellow/.style={tensor, fill=orange!80},
  avg/.style={tensor, fill=green!70!blue},
  wire/.style={thick},
}
    \node (x) at (0,2.3) {};
    \node (a1) at (-2,0.7) {};
    \node[avg] (a2) at (0,0.7) {};
    \node (a3) at (2,0.7) {};
    \node (y1) at (-2,-1) {};
    \node[yellow] (y2) at (0,-1) {};
    \node (y3) at (2,-1) {};
    \draw[wire] (a1) -- (a2);
    \draw[wire] (a2) -- (a3);
    \draw[wire] (y1) -- (y2);
    \draw[wire] (y2) -- (y3);
    \draw[wire] (a2) -- (y2);
    \draw[wire] (a2) -- (x);
    \node at (0.6,1.9) {\scriptsize $s\!+\!\ell$};
    \node at (0.45,-0.2) {\scriptsize $n_q$};
    \node at (-1,1.1) {\scriptsize $r_q^A$};
    \node at (1.3,1.1) {\scriptsize $r_{q+1}^A$};
    \node at (-1,-1.5) {\scriptsize $r_q^v$};
    \node at (1.3,-1.5) {\scriptsize $r_{q+1}^v$};
    \draw[draw=black, fill=none, rounded corners=2pt, thick] (-0.5,1.2) rectangle (0.5,-1.45);
\end{tikzpicture}\hspace{20mm}\begin{tikzpicture}[baseline=-0.5ex, scale=0.48]
\tikzset{
  tensor/.style={circle, draw, thick, minimum size=3mm},
  green/.style={tensor, fill=green!60},
  blue/.style={tensor, fill=cyan!60},
  yellow/.style={tensor, fill=orange!80},
  avg/.style={tensor, fill=green!70!blue},
  wire/.style={thick},
}
    \node (x) at (0,2.3) {};
    \node (a1) at (-2,0.7) {};
    \node[avg] (a2) at (0,0.7) {};
    \node (a3) at (2,0.7) {};
    \node (y1) at (-2,-1) {};
    \node[yellow] (y2) at (0,-1) {};
    \node (y3) at (2,-1) {};
    \draw[wire] (a1) -- (a2);
    \draw[wire] (a2) -- (a3);
    \draw[wire] (y1) -- (y2);
    \draw[wire] (y2) -- (y3);
    \draw[wire] (a2) -- (y2);
    \draw[wire] (a2) -- (x);
    \node at (0.55,1.8) {\scriptsize $m_q$};
    \node at (0.45,-0.2) {\scriptsize $n_q$};
    \node at (-1,1.1) {\scriptsize $r_q^A$};
    \node at (1.3,1.1) {\scriptsize $r_{q+1}^A$};
    \node at (-1,-1.5) {\scriptsize $r_q^v$};
    \node at (1.3,-1.5) {\scriptsize $r_{q+1}^v$};
    \draw[draw=black, fill=none, rounded corners=2pt, thick] (-0.5,1.2) rectangle (0.5,-1.45);
\end{tikzpicture}
\caption{Unfavorable interior contractions in the pre-sketched setting, which couple the TT ranks of $A$ and $v$ simultaneously and are thus quadratic in both. Left: The naive contraction arising for $\Omega_k$ (and for the sites $q \neq k$ of $\Psi_k$). Right: The $k$-th site contraction for $\Psi_k$, in which the row mode $m_k$ remains free. The rectangular boxes indicate the indices contracted at a single site.}
\label{fig:node_contraction_cost}
\end{figure}

\begin{figure}[ht]
\centering
\begin{tikzpicture}[baseline=-0.5ex, scale=0.48]
\tikzset{
  tensor/.style={circle, draw, thick, minimum size=3mm},
  green/.style={tensor, fill=green!60},
  blue/.style={tensor, fill=cyan!60},
  yellow/.style={tensor, fill=orange!80},
  avg/.style={tensor, fill=green!70!blue},
  wire/.style={thick},
}
    \node (t1) at (-2, 2.8) {};
    \node (t2) at (0, 2.8) {};
    \node (v1) at (-2, 1.7) {};
    \node (v2) at (0, 1.7) {};
    \node (k1) at (-2.3, 2.55) {};
    \node (k2) at (0.3, 2.55) {};
    \draw[wire] (-1, 2.55) -- (-1, 3);
    \draw[wire] (-1, 3) -- (3, 3);
    \draw[wire] (t1) -- (v1);
    \draw[wire] (t2) -- (v2);
    \draw[wire] (k1) -- (k2);
    \draw[wire] (a3) -- (3,0.7);
    \draw[wire] (y3) -- (3,-1);
    \node[circle, draw, fill=white, inner sep=1.5pt] at (-1,2.55) {};
    \node (x1) at (-2,2.3) {};
    \node (x2) at (0,2.3) {};
    \node[avg] (a1) at (-2,0.7) {};
    \node[avg] (a2) at (0,0.7) {};
    \node[avg] (a3) at (2,0.7) {};
    \node[yellow] (y1) at (-2,-1) {};
    \node[yellow] (y2) at (0,-1) {};
    \node[yellow] (y3) at (2,-1) {};
    \draw[wire] (a3) -- (2, 3);
    \draw[wire] (a3) -- (y3);
    \draw[wire] (a1) -- (a2);
    \draw[wire] (a2) -- (a3);
    \draw[wire] (y1) -- (y2);
    \draw[wire] (y2) -- (y3);
    \draw[wire] (a1) -- (y1);
    \draw[wire] (a2) -- (y2);
    \draw[wire] (a1) -- (x1);
    \draw[wire] (a2) -- (x2);
    \node[circle, draw, fill=white, inner sep=1.5pt] at (2,3) {};
    \node at (-1.4,1.9) {\scriptsize $s\!+\!\ell$};
    \node at (0.6,1.9) {\scriptsize $s\!+\!\ell$};
    \node at (-1.55,-0.2) {\scriptsize $n_1$};
    \node at (0.45,-0.2) {\scriptsize $n_2$};
      \node at (-1,1.1) {\scriptsize $r_1^A$};
    \node at (1.1,1.1) {\scriptsize $r_{2}^A$};
    \node at (-1,-1.5) {\scriptsize $r_1^v$};
    \node at (1.1,-1.5) {\scriptsize $r_{2}^v$};
      \draw[draw=black, fill=none, rounded corners=2pt, thick] (-2.5,1.2) rectangle (-1.5,-1.45);
\end{tikzpicture}\hspace{6mm}
\begin{tikzpicture}[baseline=-0.5ex, scale=0.48]
\tikzset{
  tensor/.style={circle, draw, thick, minimum size=3mm},
  green/.style={tensor, fill=green!60},
  blue/.style={tensor, fill=cyan!60},
  yellow/.style={tensor, fill=orange!80},
  grey/.style={tensor, fill=black!30},
  avg/.style={tensor, fill=green!70!blue},
  wire/.style={thick},
}
    \node (t2) at (0, 2.8) {};
    \node (v2) at (0, 1.7) {};
    \node (k2) at (0.3, 2.55) {};
    \draw[wire] (-1, 3) -- (3, 3);
    \draw[wire] (-1, 2.55) -- (-1, 3);
    \draw[wire] (t1) -- (v1);
    \draw[wire] (t2) -- (v2);
    \draw[wire] (k1) -- (k2);
    \node[circle, draw, fill=white, inner sep=1.5pt] at (-1,2.55) {};
    \node (x1) at (-2,2.3) {};
    \node (x2) at (0,2.3) {};
    \node[grey] (a1) at (-2,0.7) {};
    \node[avg] (a2) at (0,0.7) {};
    \node[avg] (a3) at (2,0.7) {};
    \node[yellow] (y2) at (0,-1) {};
    \node[yellow] (y3) at (2,-1) {};
    \draw[wire] (a3) -- (y3);
    \draw[wire] (a3) -- (2, 3);
    \draw[wire] (a1) -- (a2);
    \draw[wire] (a2) -- (a3);
    \draw[wire] (a1) -- (y2);
    \draw[wire] (y2) -- (y3);
    \draw[wire] (a2) -- (y2);
    \draw[wire] (a1) -- (x1);
    \draw[wire] (a2) -- (x2);
    \draw[wire] (a3) -- (3,0.7);
    \draw[wire] (y3) -- (3,-1);
    \node[circle, draw, fill=white, inner sep=1.5pt] at (2,3) {};
    \node at (-1.4,1.9) {\scriptsize $s\!+\!\ell$};
    \node at (0.6,1.9) {\scriptsize $s\!+\!\ell$};
    \node at (-1.05,-0.25) {\scriptsize $r_1^v$};
    \node at (0.45,-0.2) {\scriptsize $n_2$};
      \node at (-1,1.1) {\scriptsize $r_1^A$};
    \node at (1.1,1.1) {\scriptsize $r_{2}^A$};
    \node at (1.1,-1.5) {\scriptsize $r_{2}^v$};
      \draw[draw=black, fill=none, rounded corners=2pt, thick, 
      rotate around={50:(-1, -0.125)}] 
      (-1.55, 1.7) rectangle (-0.45, -1.9);
\end{tikzpicture}\hspace{6mm}
\begin{tikzpicture}[baseline=-0.5ex, scale=0.48]
\tikzset{
  tensor/.style={circle, draw, thick, minimum size=3mm},
  green/.style={tensor, fill=green!60},
  blue/.style={tensor, fill=cyan!60},
  yellow/.style={tensor, fill=orange!80},
  grey/.style={tensor, fill=black!30},
  avg/.style={tensor, fill=green!70!blue},
  wire/.style={thick},
}

    \node (t1) at (-2, 2.8) {};
    \node (t2) at (0, 2.8) {};
    \node (v1) at (-2, 1.7) {};
    \node (v2) at (0, 1.7) {};
    \node (k1) at (-2.3, 2.55) {};
    \node (k2) at (0.3, 2.55) {};
    \draw[wire] (-1, 3) -- (3, 3);
    \draw[wire] (-1, 2.55) -- (-1, 3);
    \draw[wire] (t1) -- (v1);
    \draw[wire] (t2) -- (v2);
    \draw[wire] (k1) -- (k2);
    \node[circle, draw, fill=white, inner sep=1.5pt] at (-1,2.55) {};
    \node (x1) at (-2,2.05) {};
    \node (x2) at (0,2.3) {};
    \node (x3) at (-2,-1) {};
    \node (x4) at (-0.3,-1) {};
    \node[avg] (a2) at (0,0.7) {};
    \node[avg] (a3) at (2,0.7) {};
    \node[grey] (y2) at (0,-1) {};
    \node[yellow] (y3) at (2,-1) {};
    \draw[wire] (a3) -- (y3);
    \draw[wire] (a3) -- (2, 3);
    \draw[wire] (a2) -- (a3);
    \draw[wire] (y2) -- (y3);
    \draw[wire] (a2) -- (y2);
    \draw[wire] (a2) -- (x2);
    \draw[wire] (a3) -- (3,0.7);
    \draw[wire] (y3) -- (3,-1);
    \draw[wire] plot [smooth, tension=0] coordinates {(x1) (x3) (x4)};
    \node[circle, draw, fill=white, inner sep=1.5pt] at (2,3) {};
    \node at (-1.4,1.9) {\scriptsize $s\!+\!\ell$};
    \node at (0.6,1.9) {\scriptsize $s\!+\!\ell$};
    \node at (0.8,-0.2) {\scriptsize $n_2 r_1^A$};
    \node at (1.1,1.1) {\scriptsize $r_{2}^A$};
    \node at (1.1,-1.5) {\scriptsize $r_{2}^v$};
      \draw[draw=black, fill=none, rounded corners=2pt, thick] 
      (-2.3, 2.85) rectangle (0.5, -1.45);
\end{tikzpicture}\hspace{6mm}\begin{tikzpicture}[baseline=-0.5ex, scale=0.48]
\tikzset{
  tensor/.style={circle, draw, thick, minimum size=3mm},
  green/.style={tensor, fill=green!60},
  blue/.style={tensor, fill=cyan!60},
  yellow/.style={tensor, fill=orange!80},
  grey/.style={tensor, fill=gray!60},
  avg/.style={tensor, fill=green!70!blue},
  wire/.style={thick},
}
    \node (t2) at (0, 2.8) {};
    \node (v2) at (0, 1.7) {};
    \node (k2) at (0.3, 2.55) {};
    \draw[wire] (0, 3) -- (3, 3);
    \node (x2) at (0,2.3) {};
    \node (x3) at (-2,-1) {};
    \node (x4) at (-0.3,-1) {};
    \node[grey] (a2) at (0,0.7) {};
    \node[avg] (a3) at (2,0.7) {};
    \node[yellow] (y3) at (2,-1) {};
    \draw[wire] (a3) -- (y3);
    \draw[wire] (a3) -- (2, 3);
    \draw[wire] (a2) -- (a3);
    \draw[wire] (a2) -- (y3);
    \draw[wire] (a2) -- (0,3);
    \draw[wire] (a3) -- (3,0.7);
    \draw[wire] (y3) -- (3,-1);
    \node[circle, draw, fill=white, inner sep=1.5pt] at (2,3) {};
    \node at (0.6,1.9) {\scriptsize $s\!+\!\ell$};
    \node at (1.1,1.1) {\scriptsize $r_{2}^A$};
        \node at (0.95,-0.25) {\scriptsize $r_2^v$};
    \draw[draw=black, fill=none, rounded corners=2pt, thick, 
      rotate around={50:(1, -0.125)}] 
      (0.45, 1.7) rectangle (1.55, -1.9);
\end{tikzpicture}
\caption{Efficient left-to-right contraction sequence. Contractions are performed site by site so that only one TT rank (of $A$ or $v$) contributes quadratically at each step, avoiding joint rank coupling.}
\label{fig:left-to-right_contraction}
\end{figure}
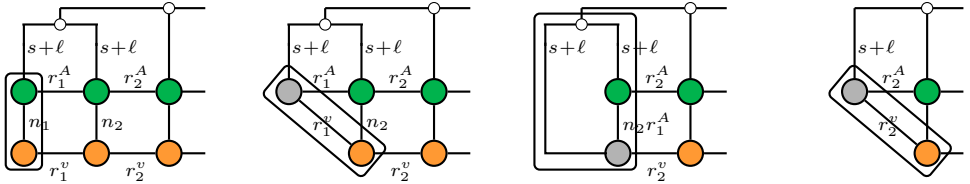

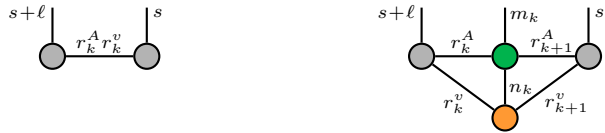
\begin{figure}[ht]
\centering
\begin{tikzpicture}[baseline=-0.5ex, scale=0.48]
\tikzset{
  tensor/.style={circle, draw, thick, minimum size=3mm},
  green/.style={tensor, fill=green!60},
  grey/.style={tensor, fill=gray!60},
  blue/.style={tensor, fill=cyan!60},
  yellow/.style={tensor, fill=orange!80},
  avg/.style={tensor, fill=green!70!blue},
  wire/.style={thick},
}
    \node (x1) at (-2.3,2.3) {};
    \node (x2) at (0.3,2.3) {};
    \node[grey] (a1) at (-2.3,0.7) {};
    \node[grey] (a2) at (0.3,0.7) {};
    \node (a3) at (2,0.7) {};
    \draw[wire] (a1) -- (a2);
    \draw[wire] (a1) -- (x1);
    \draw[wire] (a2) -- (x2);
    \node at (0.6,1.9) {\scriptsize $s$};
    \node at (-3,1.9) {\scriptsize $s\!+\!\ell$};
      \node at (-1,1.1) {\scriptsize $r_k^A r_k^v$};
\end{tikzpicture}\hspace{20mm}\begin{tikzpicture}[baseline=-0.5ex, scale=0.48]
\tikzset{
  tensor/.style={circle, draw, thick, minimum size=3mm},
  green/.style={tensor, fill=green!60},
  blue/.style={tensor, fill=cyan!60},
  yellow/.style={tensor, fill=orange!80},
  avg/.style={tensor, fill=green!70!blue},
  wire/.style={thick},
  grey/.style={tensor, fill=gray!60},
}
    \node (x1) at (-2.3,2.3) {};
    \node (x2) at (0,2.3) {};
    \node (x3) at (2.3,2.3) {};
    \node[grey] (a1) at (-2.3,0.7) {};
    \node[avg] (a2) at (0,0.7) {};
    \node[grey] (a3) at (2.3,0.7) {};
    \node (y1) at (-2,-1) {};
    \node[yellow] (y2) at (0,-1) {};
    \draw[wire] (a1) -- (a2);
    \draw[wire] (a1) -- (y2);
    \draw[wire] (a3) -- (y2);
    \draw[wire] (a2) -- (a3);
    \draw[wire] (a2) -- (y2);
    \draw[wire] (a1) -- (x1);
    \draw[wire] (a2) -- (x2);
    \draw[wire] (a3) -- (x3);
        \node at (-3,1.9) {\scriptsize $s\!+\!\ell$};
            \node at (2.6,1.9) {\scriptsize $s$};
    \node at (0.55,1.8) {\scriptsize $m_k$};
    \node at (0.45,-0.2) {\scriptsize $n_k$};
      \node at (-1.15,1.1) {\scriptsize $r_k^A$};
    \node at (1.3,1.1) {\scriptsize $r_{k+1}^A$};
    \node at (-1.4,-0.6) {\scriptsize $r_k^v$};
    \node at (1.7,-0.6) {\scriptsize $r_{k+1}^v$};
\end{tikzpicture}
\caption{Illustration of the final step in the compression of $\Omega_k$ (left) and $\Psi_k$ (right).}
\label{fig:final_contraction}
\end{figure}

Sweeping sequentially from the peripheral cores is beneficial across all three routes: it allows intermediate contractions to be recycled and keeps the memory footprint strictly bounded. For the fused and pre-sketched routes, sweeping also brings a major efficiency gain by avoiding isolated contractions at interior sites, which would otherwise couple the TT ranks of $A$ and $v$ simultaneously and yield a cost that is quadratic in both. We illustrate the sweeping mechanism in the context of the pre-sketched approach. \Cref{fig:node_contraction_cost} shows the two costly scenarios arising from naive interior contractions, while \cref{fig:left-to-right_contraction} demonstrates how the left-to-right sweep circumvents them by processing one rank at a time. Once both sweeps are complete, $\Omega_k$ and $\Psi_k$ are assembled via the merging contractions depicted in \cref{fig:final_contraction}.


The two families are not symmetric, and this asymmetry is what limits how much
the pre-sketch can help. The matrix $\Omega_k$ is sketched on both sides, so it can
be assembled from the sketched cores alone. The matrix $\Psi_k$ instead leaves the
$k$-th row mode free, because that mode carries the $k$-th core of the
approximant and no sketch is applied to it; its assembly therefore requires
$G_k$ in its original form. Every contraction of the algorithm except this one
sees the operator only through $\widetilde G_k$, and costs $\mathcal{O}(n_k)$ per
unit of the surrounding ranks instead of $\mathcal{O}(m_kn_k)$; the assembly of
$\Psi_k$ is the single place where the full $m_k\times n_k$ block is still
traversed.

\subsection{Cost}
\label{sec:cost}
To compare the three orderings, we assume uniform mode sizes $m_k=m$, $n_k=n$,
uniform TT ranks $r^A_k=r^A$, $r^v_k=r^v$, target rank $s$, and
oversampling parameter $\ell$. We recall that we focus on the regime of large
mode sizes~$m$ and~$n$, in which the compression phase dominates throughout:
recovering the TT-cores afterwards requires only the singular value
decompositions of the $d-1$ matrices~$\Omega_k$, of size
$(s+\ell)\times s$, and the application of the resulting pseudoinverses to
the $\Psi_k$; neither operation involves an $m\times n$ block, so both are
lower-order in the mode sizes.

The three orderings perform several contractions, but only some of them touch
the $m\times n$ blocks of the cores of $A$, which dominate the cost when
$m$ and $n$ are large. The remaining terms are lower-order in this regime
and do not affect the comparison. We therefore report, for each ordering,
only the terms containing $m\,n$. Each dominant cost is a multiple of
$d\,m\,n\,(r^A)^2$, and what changes from one ordering to another is that
multiple.

\paragraph{The sequential ordering} This contracts each core of $A$ with the
corresponding core of $v$ to form the cores of $Av$, of which the TT ranks are
$r^A r^v$, and compresses them one at a time. This formation constitutes the dominant
cost, which is
\[
  \mathcal{O}\!\left(d\,m\,n\,(r^A)^2(r^v)^2\right),
\]
so the $m\times n$ blocks are traversed with weight $(r^v)^2$. Nothing here depends on the target rank: the same work is required whether a low-rank or high-rank approximation is sought. 
When memory is not a limiting factor, the cores of $Av$ can be stored and reused in the two sweeps, so that they need not be formed twice. This halves the dominant cost associated with their formation. In the experiments reported below, however, we use the memory-saving variant, which avoids storing the cores of $Av$.
\paragraph{The fused ordering}
First, each core of $v$ is contracted with the partial contraction arriving from the previous site. The result is then contracted with the corresponding core of $A$ and, finally, with the sketch. Only the second contraction involves the $m\times n$ block. Depending on whether the sweep proceeds from the right or from the left, its cost is
\[
\mathcal{O}\!\left(d\,m\,n\,(r^A)^2\,r^v s\right)
\quad\text{or}\quad
\mathcal{O}\!\left(d\,m\,n\,(r^A)^2\,r^v(s+\ell)\right),
\]
respectively. Note that assembling the $\Psi_k$ requires no additional pass over $A$, since the preceding contractions can be recycled. The combined cost of the two sweeps therefore has a weight of $r^v(2s+\ell)$.

\paragraph{The pre-sketched ordering} This contracts the sketching matrices into
the cores of $A$ before the operator is applied, so that the sweeps no longer
see a row mode at all. A single contraction per product still contains $m\,n$,
the assembly of the $\Psi_k$,
\[
  \mathcal{O}\!\left(d\,m\,n\,(r^A)^2\,s\,(s+\ell)\right),
\]
for a weight of $s\,(s+\ell)$. The absorption of the sketch is itself a traversal
of the $m\times n$ blocks, of weight $s+\ell$, but it does not involve $v$ and
is therefore carried out once per operator rather than once per product; even when it
is carried out every time the cost stays a factor $s$ below the assembly of the $\Psi_k$.

\begin{table}[ht]
\centering
\caption{Weight multiplying $d\,m\,n\,(r^A)^2$ in the cost of traversing the
$m\times n$ blocks of the cores of $A$.
The second column counts the work repeated for every product; the third counts
the work that depends on $A$ alone and is therefore performed once, however
many vectors the operator is applied to.}
\label{tab:weights}
\begin{tabular}{lcc}
\hline
ordering & per product & once per operator \\
\hline
sequential     & $2(r^v)^2$                   & --- \\
fused          & $r^v\left(2s+\ell\right)$ & --- \\
pre-sketched   & $s\left(s+\ell\right)$      & $s+\ell$ \\
\hline
\end{tabular}
\end{table}

\cref{tab:weights} summarizes the comparison:
 for large mode sizes $m$ and $n$, the choice between the different orderings is largely determined by the relative sizes of $s$ and~$r^v$: the ratio between the sequential and pre-sketched costs is of order $\left(r^v/s\right)^2$. Thus, pre-sketching is particularly advantageous when the desired approximation rank $s$ is small compared to the input TT rank $r^v$, the sequential ordering becomes preferable as $s$ approaches or exceeds $r^v$, and the fused ordering interpolates between the two extremes.
 However, if the TT ranks become large compared to the mode sizes $m$ and~$n$, other terms in the cost can become dominant, and the estimates in \cref{tab:weights} no longer capture the full advantage of the pre-sketched ordering.

For problems with small mode sizes, such as those arising in quantum chemistry and many-body quantum dynamics, the computational advantage of pre-sketching is reduced. The sequential and fused variants, however, remain effective streamable methods, combining low memory requirements with the theoretical guarantees of randomized sketching.

\subsection{Memory}
\label{sec:memory}

Which objects have to be stored depends on the setting: the operator may already
be in memory or have to be built, the argument may be stored or generated on the
fly, and both may carry structure. We limit the comparison to the dense case.
\begin{table}[ht]
\centering
\caption{Storage, in numbers, of the objects involved in one product. The
sketched operator and the sketches are those of the pre-sketched ordering; the
cores of $Av$ are what a route that forms the product must hold.}
\label{tab:memory}
\begin{tabular}{ll}
\hline
object & storage \\
\hline
operator $A$                      & $\left(2r^A + (d-2)(r^A)^2\right)m\,n$ \\
argument $v$                      & $\left(2r^v + (d-2)(r^v)^2\right)n$ \\
one interior core of $Av$         & $m\,(r^A r^v)^2$ \\
all the cores of $Av$             & $\left(2r^Ar^v + (d-2)(r^Ar^v)^2\right)m$ \\
sketches of the operator                 & $\left(2r^A + (d-2)(r^A)^2\right)n\,(s+\ell)$ \\
sketches $\Omega_k$, $\Psi_k$     & $(d-1)(s+\ell)s + \left(2s+\ell + (d-2)(s+\ell)s\right)m$ \\
\hline
\end{tabular}
\end{table}
Table~\ref{tab:memory} summarizes the memory costs. The sketched operator
requires a fraction $(s+\ell)/m$ of the storage of $A$, and at each interior
site the sketch $\Psi_k$, which takes the place of the $k$-th core of $Av$,
is smaller than that core by a factor $(s+\ell)\,s/(r^A r^v)^2$. More importantly, the sequential
ordering must still materialize at least one interior core of $Av$, requiring
$m\,(r^A r^v)^2$ numbers. This can make it infeasible when $r^A r^v$ is large.

The fused and pre-sketched orderings avoid forming cores of $Av$. The former
requires an intermediate of size $r^A n \, r^v s$, while the latter stores the
sketched operator. The pre-sketch is reusable across multiple products with
the same operator, whereas the fused ordering requires no additional
persistent storage. Both intermediates can also be processed in blocks of
sketch columns to reduce the working memory.

For structured operators, however, pre-sketching may be less attractive:
sketching destroys sparsity and produces dense cores. 

\subsection{Comparison with SRC}
\label{sec:src}

In the regime considered here, the straightforward implementation of SRC described in \cite{camano2026successive}, which we refer to as \emph{standard SRC}, is more expensive than our pre-sketched approach. The approach in \cite{camano2026successive} is natural in the regime considered therein, but when the mode sizes are large and $s\ll r^v$, pre-sketched approaches become advantageous. In fact, the first sweep of SRC can be performed in the same way as the first sweep of STTA, by absorbing the sketch into the operator cores in advance. A similar contraction ordering can be used in the second sweep, by contracting the orthonormal factors computed by the local QR factorizations into the operator cores before the remaining contractions. These factors, however, depend on the input and must therefore be recomputed for every matvec. We refer to this reorganized version as \emph{pre-sketched SRC}. We include both implementations of SRC in the comparisons below.

\subsection{Repeated application}
\label{sec:krylov}

Iterative methods provide a natural setting in which the proposed operator recompression can be applied repeatedly. This is particularly relevant in the TT format, where successive operator applications may change or increase the ranks of the iterates.

Without truncation, repeated operator applications cause ranks to grow steadily. However, pre-sketching provides a clear computational advantage only under aggressive truncation; specifically, when the target compression rank is significantly smaller than the rank of the input vector ($s \ll r^v$).

We therefore explored whether this trade-off could be exploited by allowing TT ranks to grow between truncations. In particular, we tested:
\begin{itemize}
\item allowing ranks to grow up to a threshold $r_{\max}$ over several iterations before applying a sharp truncation to $r_{\min}$;
\item applying rank-inflating operations, such as accurate preconditioning ($PAv$ or $APv$) or non-linear operations, followed by aggressive truncation to $r_{\min}$.
\end{itemize}
Experimentally, these strategies proved advantageous only in narrow, somewhat artificial scenarios. In most of our test cases, maintaining a fixed rank $r_{\min}$ throughout, with truncation after every iteration, achieved comparable or superior convergence at a lower overall computational cost.

\section{Numerical experiments}
\label{sec:experiments}

All timings were obtained in MATLAB R2026a on an Intel Core Ultra 5 235U
processor with 12 physical cores, using the Intel oneAPI Math Kernel Library
2025.0.1, as well as the TT-Toolbox \cite{oseledets2011tensor}.  Each reported timing
is the fastest of five runs.  The code and data
generating each figure are available at
\url{https://github.com/alb95/fast-tt-matvec}.
\subsection{Synthetic experiments}
\label{sec:synthetic_exp}

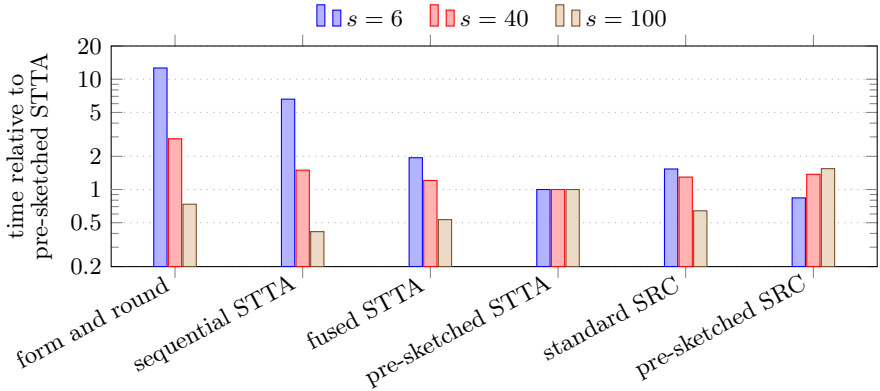
\begin{figure}[tbp]
  \centering
  \begin{tikzpicture}
    \begin{axis}[
      width=0.9\linewidth, height=4.5cm,
      ybar=0.5pt, bar width=5pt,
      ymode=log, log origin=infty,
      ymin=0.2, ymax=20,
      enlarge x limits=0.10,
      symbolic x coords={formround,sequential,fused,presketched,src,srcabsorbed},
      xtick=data,
      xticklabels={{form and round},{sequential STTA},{fused STTA},{pre-sketched STTA},%
                   {standard SRC},{pre-sketched SRC}},
      x tick label style={rotate=28, anchor=east, font=\small},
      ylabel={\shortstack{time relative to\\pre-sketched STTA}},
      ylabel style={font=\small},
      ytick={0.2,0.5,2,5,10,20},
      yticklabels={$0.2$,$0.5$,$2$,$5$,$10$,$20$},
      minor ytick={0.3,0.4,0.6,0.7,0.8,0.9,3,4,6,7,8,9,30},
      yticklabel style={font=\small},
      ymajorgrids, major grid style={dotted, black!30},
      extra y ticks={1}, extra y tick labels={$1$},
      extra y tick style={grid=major, grid style={solid, black!55}},
      legend style={font=\small, at={(0.5,1.02)}, anchor=south,
                    legend columns=3, draw=none, fill=none,
                    /tikz/every even column/.append style={column sep=0.8em}},
      legend cell align=left,
      ]
      \addplot table [x=route, y=s6, col sep=comma] {exp1_target_rank_plot.csv};
      \addlegendentry{$s=6$}
      \addplot table [x=route, y=s40, col sep=comma] {exp1_target_rank_plot.csv};
      \addlegendentry{$s=40$}
      \addplot table [x=route, y=s100, col sep=comma] {exp1_target_rank_plot.csv};
      \addlegendentry{$s=100$}
    \end{axis}
  \end{tikzpicture}

  \caption[Running time of each route at three target ranks]{%
    Running time of each route divided by the running time of
    pre-sketched STTA, for three target ranks, $s=6,40,100$, and everything else fixed:
    $d=4$, $m=n=512$, $r^A=4$, $r^v=40$, and $\ell=4$.  
  }
  \label{fig:exp-target-rank}
\end{figure}

The first experiment compares the computational cost of six approaches for
computing a low-rank approximation of $Av$: forming $Av$ explicitly and
rounding it, the sequential, fused, and pre-sketched STTA orderings, and the
standard and pre-sketched SRC approaches described in
\cref{sec:src}. The operator $A$, the vector $v$, both with dense cores,
and the random sketching matrices are drawn randomly and kept fixed; only the
target rank $s$ changes, with the sketching matrices nested as $s$ increases.
Thus, within each group of bars in \cref{fig:exp-target-rank}, all quantities
other than $s$ are identical.

We choose $m=n=512$ and $r^v=40$ to obtain a problem that is sufficiently
large for the costs associated with the cores of $A$ to dominate, while
remaining feasible for all six approaches and allowing target ranks both well
below and above $r^v$. At these sizes, however, secondary costs can still
affect the timings, so larger problems would provide a more reliable
assessment of the predicted scaling. Although the reported timings are
averaged over multiple runs, their quantitative comparison is also affected by
how efficiently MATLAB handles the different matrix operations at the
resulting sizes.

Overall, the results largely follow the scaling predicted by
\cref{tab:weights}: the pre-sketched methods are fastest for small $s$, while
the sequential ordering becomes the fastest for large $s$, with the transition
occurring around $s\approx r^v$.

\pgfplotsset{
  ttstyle/.style={
    width=0.47\linewidth, height=0.38\linewidth,
    tick label style={font=\scriptsize},
    label style={font=\small},
    legend style={font=\scriptsize, draw=none, fill=none, at={(0.02,0.98)},
                  anchor=north west, cells={anchor=west}},
    grid=both, grid style={line width=0.2pt, draw=black!12},
    every axis plot/.append style={line width=0.8pt, mark size=1.6pt},
  },
  form/.style ={black!55, mark=square*,    mark options={fill=black!55}},
  seqs/.style ={black!55, mark=square,     dashed},
  fuse/.style ={blue!60!black,  mark=triangle*},
  pres/.style ={red!70!black,   mark=*},
  srcs/.style ={green!45!black, mark=diamond*},
}
\subsection{A near-resonant solve}
\label{sec:resonance}

We have seen in the previous experiments that the pre-sketched algorithm is
cheaper than the other alternatives when $s\ll r^v$.
Nevertheless, it is not obvious for which problems this situation arises.
A particularly simple example is the solution of a nearly singular linear
system, written as the application of its inverse: $A$ is then an operator
that strongly amplifies one direction, and that direction may have only a
small component in $v$.
Suppose that $u$ is such a direction and write
$v = c u + v_\perp$, $|c|\ll 1$, with $v_\perp$ orthogonal to $u$. If $A$
amplifies $u$ by a factor much larger than the other directions, then
$Av \approx c\,Au + A v_\perp$ can be dominated by the first term even though
$c$ is small. This creates an important mismatch: the small component $cu$ can
be essential for computing $Av$, while the result $Av$ itself can be well
approximated by a low-rank object.

This is precisely where compressing $v$ before applying $A$ can fail. If a
rank-$s$ approximation of $v$ discards the small component in the direction
$u$, it may remove the component that determines the solution. The alternative
of first forming $Av$ and then compressing it is mathematically correct, but
can be prohibitively expensive when $A$ and $v$ are represented in tensor
format. Our primitive targets the intermediate possibility: apply $A$ and
compress the result at the same time, without explicitly forming $Av$.

A near-resonant system provides a concrete and physically intuitive example
of this mechanism. Consider a symmetric operator $L$ with eigenpairs
$(\lambda_j,\varphi_j)$, and let $A=(L-\sigma I)^{-1}$ be its resolvent at a
shift $\sigma$. If $\sigma$ is close to one eigenvalue $\lambda^*$, then
the corresponding eigenvector $\varphi^*$ is strongly amplified, with
$A\varphi^*=\frac{1}{\lambda^*-\sigma}\,\varphi^*$. Thus, if $v$ has only a small
component $c=\frac{\langle \varphi^*,v\rangle}{\|\varphi^*\|^2}$
in the resonant direction $\varphi^*$, the corresponding component of $Av$
is multiplied by $(\lambda^*-\sigma)^{-1}$. By choosing $\sigma$ sufficiently
close to $\lambda^*$, this small component can dominate the whole solution.
The terminology comes from resonators: $\lambda^*$ represents a natural
mode of the system, while choosing $\sigma$ close to $\lambda^*$ corresponds
to driving the system close to that mode's natural frequency. Here, however,
the physical interpretation is secondary; what matters for us is the algebraic
mechanism: a direction that is small in the input becomes dominant after
applying~$A$.

We construct such an example in TT format. Let
$\Delta_n=h^{-2}\operatorname{tridiag}(-1,2,-1)$, with $h = \frac{1}{n+1}$ the mesh parameter, be the discretization of
$-\partial_{xx}$ with Dirichlet conditions on $n$ interior points, which is
positive definite with $\Delta_n=S\diag(\mu_1,\dots,\mu_n)S^\top$, and let $L$
be the Kronecker sum of $d$ copies of it. Its eigenvectors are separable and its
eigenvalues are the sums $\mu_{j_1}+\cdots+\mu_{j_d}$, so $L$ is again positive
definite. For $\sigma$ inside the spectrum the shifted matrix $L-\sigma I$ is
symmetric indefinite and nonsingular; writing $\delta$ for its eigenvalue of
smallest modulus, $A=(L-\sigma I)^{-1}$ is symmetric with
$\|A\|_2=\delta^{-1}$.

The resolvent itself is available in tensor format, and for the same reason
that $L$ is. Being a Kronecker sum, $L$ is diagonalized by the orthogonal
factor $S\otimes\cdots\otimes S$, so that
\begin{equation*}
  A=(S\otimes\cdots\otimes S)\,\diag(\mathcal F)\,
     (S\otimes\cdots\otimes S)^\top,
  \qquad
  \mathcal F_{j_1\dots j_d}=\frac{1}{\mu_{j_1}+\cdots+\mu_{j_d}-\sigma}.
\end{equation*}
The outer factors are rank-one TT operators, so the TT ranks of $A$ are those
of the $d$-dimensional array $\mathcal F$, which we compute by TT-SVD and
truncate at tolerance $10^{-8}$, giving $r^A\le16$.  Those ranks stay small
because $\mathcal F$ is a smooth function of the sum
$\mu_{j_1}+\cdots+\mu_{j_d}$ except at the resonance, where it has a single
spike.  At this cap the truncation of $\mathcal F$ contributes about
$10^{-8}$ in relative Frobenius norm, well below the errors reported below,
and the same $A$ is given to every route, so it is not part of the
comparison.

Products of this kind occur wherever a resolvent is applied at a shift close
to the spectrum: in the evaluation of frequency responses, and in the
correction-vector step by which spectral functions are computed in tensor
formats \cite{jeckelmann2008density,nocera2016spectral}, as well as in the
shifted systems of parametrized problems \cite{kressner2010krylov}.  What
makes $w=Av$ of low rank here is that the eigenvectors of $L$ are separable:
for $\lambda^*=\mu_{j_1}+\mu_{j_2}+\mu_{j_3}$ with $j_1$, $j_2$, $j_3$ distinct the
resonant eigenspace has multiplicity $d!=6$, is spanned by the permutations of
$\varphi_{j_1}\otimes\varphi_{j_2}\otimes\varphi_{j_3}$, and has TT rank at
most $3$. Truncating $w$ at that rank therefore leaves out only its
non-resonant part, the size of which is proportional to $\delta$.

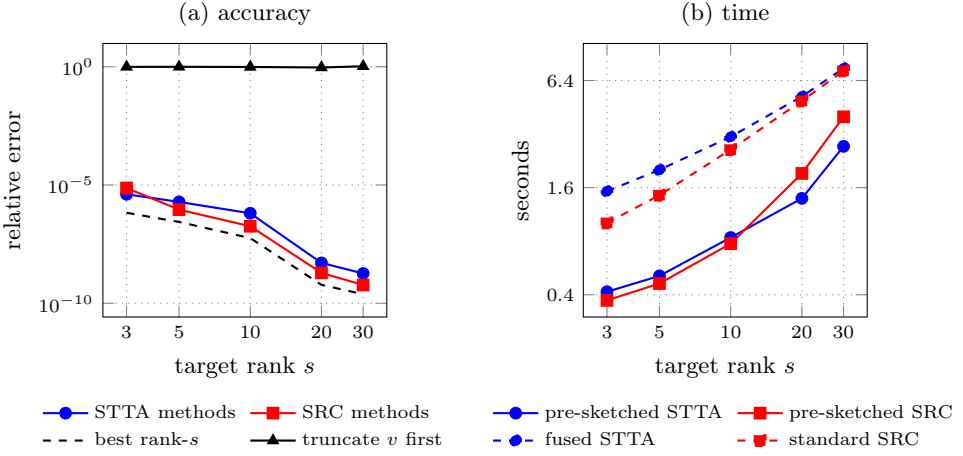
\begin{figure}[tbp]
  \centering
  \pgfplotstableread[col sep=comma]{response_time.csv}\responsedata

  \begin{tikzpicture}
    \begin{groupplot}[
      group style={
        group size=2 by 1,
        horizontal sep=2.6cm,
      },
      width=0.41\linewidth,
      height=5.2cm,
      xmode=log,
      log basis x=2,
      xtick=data,
      xticklabels from table={\responsedata}{s},
      xlabel={target rank $s$},
      xlabel style={font=\small},
      ylabel style={font=\small},
      tick label style={font=\scriptsize},
      unbounded coords=jump,
      grid=major,
      grid style={dotted, black!50},
      legend style={
        font=\scriptsize,
        draw=none,
        fill=none,
        at={(0.5,-0.28)},
        anchor=north,
        /tikz/every even column/.append style={
          column sep=0.3em
        },
      },
      legend cell align=left,
      every mark/.append style={scale=0.7},
    ]

    \nextgroupplot[
      ymode=log,
      ylabel={relative error},
      title={\small (a) accuracy},
      title style={yshift=-1mm},
      legend columns=2,
    ]

    \addplot+[thick, blue, mark=*]
      table [x=s, y=e_pre] {\responsedata};
    \addlegendentry{STTA methods}

    \addplot+[thick, red, mark=square*]
      table [x=s, y=e_src] {\responsedata};
    \addlegendentry{SRC methods}

    \addplot+[thick, black, dashed, mark=none]
      table [x=s, y=e_floor] {\responsedata};
    \addlegendentry{best rank-$s$}

    \addplot+[thick, mark=triangle*]
      table [x=s, y=e_trunc] {\responsedata};
    \addlegendentry{truncate $v$ first}

    \nextgroupplot[
      ymode=log,
      ymin=0.3,
      ylabel={seconds},
      title={\small (b) time},
      title style={yshift=-1mm},
      ytick={0.4,1.6,6.4},
      yticklabels={0.4,1.6,6.4},
      legend columns=2,
    ]

    \addplot+[thick, blue, mark=*]
      table [x=s, y=t_pre] {\responsedata};
    \addlegendentry{pre-sketched STTA}

    \addplot+[thick, red, mark=square*]
      table [x=s, y=t_src_abs] {\responsedata};
    \addlegendentry{pre-sketched SRC}

    \addplot+[thick, blue, dashed, mark=*]
      table [x=s, y=t_fused] {\responsedata};
    \addlegendentry{fused STTA}

    \addplot+[thick, red, dashed, mark=square*]
      table [x=s, y=t_src] {\responsedata};
    \addlegendentry{standard SRC}

    \end{groupplot}
  \end{tikzpicture}

  \caption[The near-resonant solve]{%
  The product $w=Av$ with $A=(L-\sigma I)^{-1}$, for $d=3$, $n=512$,
  $r^v=100$, $r^A\le16$, $\ell=4$, and $\delta=10^{-6}\eta$; the numerical
  TT rank of $w$ at $10^{-6}$ is $3$.
  (a)~Relative error against $Av$.  The orderings within each family return
  the same approximant, so STTA and SRC are one curve each. The ``best
  rank-$s$'' approach is the TT rounding of $Av$; ``truncate $v$ first'' rounds $v$
  to rank $s$ and then applies $A$ exactly.
  (b)~Running time, with every route building its own pre-sketch inside the
  timed call and the routes timed round-robin.  Forming $Av$ is not
  attempted: at TT rank $r^Ar^v=1600$ a single interior core occupies
  $10.5$~GB, which is also why the sequential ordering is absent.%
}
  \label{fig:response}
\end{figure}

Figure~\ref{fig:response} illustrates this for $d=3$, $n=512$, and a shift at
distance $\delta=10^{-6}\eta$ from $\lambda^*$, where $\eta$ is the distance
from $\lambda^*$ to the rest of the spectrum of $L$. We take a right-hand side
$v$ of TT rank $r^v=100$, drawn at random. Its component along the resonant
eigenspace is tiny, as it is for any $v$ not built around $\varphi^*$, and it
is precisely that component which the resolvent amplifies, so that everything
else is negligible in~$w$. The experiment confirms this behavior, and the
pre-sketched ordering is three to four times faster than the fused ordering, and two to three times faster than standard SRC. It also
shows that compressing $v$ instead of the product fails: rounding $v$ to rank
$s$ can remove the very component that $A$ subsequently amplifies.
\subsection{Streaming uncertainty propagation}
\label{sec:exp-stream}

The previous experiments consider arguments that are fully available in memory. We now consider a streaming setting motivated by uncertainty quantification, where an underlying deterministic field must be recovered from repeated noisy realizations. The goal is to compute $Au$ from an ensemble of $N$ realizations
$$
v_i = u + \rho \xi_i,
$$
where $u$ is a low-rank deterministic field and $\xi_i$ is a stochastic perturbation with zero mean. In particular, assuming that the perturbations cancel in the ensemble mean, we have
$$
Au = A\left(\frac{1}{N}\sum_{i=1}^N v_i\right).
$$
The factor $1/N$ is applied once, after the stream has been processed, and plays no part in what follows. The realizations arrive sequentially and must be discarded after processing, so storing the full ensemble is not an option. The challenge is therefore to accumulate their contributions while keeping the intermediate representation compact.
 
The difficulty stems from the discrepancy between the ranks of the individual realizations, the accumulated sum, and the final solution. The deterministic field has TT rank $r^u=6$ and a rapidly decaying spectrum, whereas each stochastic perturbation has TT rank $r^\xi=50$ and consists of rough separable modes with comparable energy. Each realization consequently has TT rank $r^v=r^u+r^\xi=56$. Applying $A$ directly to a realization can produce a TT representation of rank $r^A r^v$, while their accumulation can have rank as large as $N r^A r^v=22400$. In contrast, when the perturbations cancel the ensemble mean is exactly $u$, so the solution $Au$ has TT rank at most $r^A r^u=24$. Thus, the desired result has a much lower rank than a straightforward representation of the accumulated products.
 
We set $d=4$, $n=300$, $N=100$, $r^A=4$, and $\rho=20$. The stochastic modes have comparable energy, so assuming an approximately uniform energy distribution, each has norm of order $\rho/\sqrt{r^\xi}$. Since $\|u\|_F=1$ and $\rho>\sqrt{r^\xi}\approx 7$, the stochastic modes can individually be larger than the deterministic field. This makes low-rank truncation of individual realizations particularly challenging: it preferentially retains stochastic directions and may discard components of $u$.
 
\begin{figure}[tb]
\centering
\begin{tikzpicture}
\begin{groupplot}[group style={group size=2 by 2, horizontal sep=1.6cm, vertical sep=1.8cm},
    ttstyle, 
    grid=major,
grid style={dotted, black!50},
width=0.44\linewidth, height=0.34\linewidth,
    xlabel={target rank $s$}, xlabel style={font=\small},
    ylabel style={font=\small}, tick label style={font=\small},
    legend style={font=\footnotesize, draw=none, fill=none,
                  at={(1.15,-0.37)}, anchor=north, legend columns=4,
                  /tikz/every even column/.append style={column sep=0.6em}},
    legend cell align=left, every mark/.append style={scale=0.8}]

  \nextgroupplot[ymode=log,
  ylabel={relative error},
                 title={\small (a) exact cancellation, accuracy},
                 title style={yshift=-1mm}]
  \addplot[black!45, dashed, mark=none]
      table[col sep=comma, x index=0, y index=2] {stream.csv};
  \addplot[thick, blue, mark=*]
      table[col sep=comma, x index=0, y index=3] {stream.csv};
  \addplot[thick, red, mark=square*]
      table[col sep=comma, x index=0, y index=4] {stream.csv};
  \addplot[thick, black, dashdotted, mark=triangle*]
      table[col sep=comma, x index=0, y index=5] {stream.csv};

  \nextgroupplot[
  ylabel={seconds},
                 title={\small (b) exact cancellation, time},
                 title style={yshift=-1mm}]
  \addplot[thick, blue, mark=*]
      table[col sep=comma, x index=0, y index=6] {stream.csv};
  \addplot[thick, red, mark=square*]
      table[col sep=comma, x index=0, y index=7] {stream.csv};
  \addplot[thick, black, dashdotted, mark=triangle*]
      table[col sep=comma, x index=0, y index=8] {stream.csv};

  \nextgroupplot[ymode=log, ylabel={relative error},
    yminorgrids=false,
    minor y tick style={draw=none},
                 title={\small (c) residual $10^{-3}$, accuracy},
                 title style={yshift=-1mm}]
  \addplot[black!45, dashed, mark=none]
      table[col sep=comma, x index=0, y index=2] {stream_eps.csv};
  \addlegendentry{best rank-$s$}
  \addplot[thick, blue, mark=*]
      table[col sep=comma, x index=0, y index=3] {stream_eps.csv};
  \addlegendentry{STTA, pre-sketched}
  \addplot[thick, red, mark=square*]
      table[col sep=comma, x index=0, y index=4] {stream_eps.csv};
  \addlegendentry{SRC, pre-sketched}
  \addplot[thick, black, dashdotted, mark=triangle*]
      table[col sep=comma, x index=0, y index=5] {stream_eps.csv};
  \addlegendentry{compress $v_i$ first}

  \nextgroupplot[
  ylabel={seconds},
                 title={\small (d) residual $10^{-3}$, time},
                 title style={yshift=-1mm}]
  \addplot[thick, blue, mark=*]
      table[col sep=comma, x index=0, y index=6] {stream_eps.csv};
  \addplot[thick, red, mark=square*]
      table[col sep=comma, x index=0, y index=7] {stream_eps.csv};
  \addplot[thick, black, dashdotted, mark=triangle*]
      table[col sep=comma, x index=0, y index=8] {stream_eps.csv};
\end{groupplot}
\end{tikzpicture}
\caption{One pass over $N=100$ realizations at $d=4$, $n=300$, against the
  target rank $s$.  Top row: the fluctuations cancel, so the sum is $NAu$ at
  TT rank $r^Ar^u=24$ and is represented exactly there, which is the drop at
  the last point of~(a).  Bottom row: a fluctuation of relative size
  $10^{-3}$ survives, the exact answer has TT rank $r^A(r^u+r^{\xi})=224$,
  and the best rank-$s$ error itself plateaus at about $10^{-3}$ from $s=8$
  onwards.  Both pre-sketched routes compress only after the contributions have
  been combined, and track the best rank-$s$ approximation in either case.%
}
\label{fig:stream}
\end{figure}
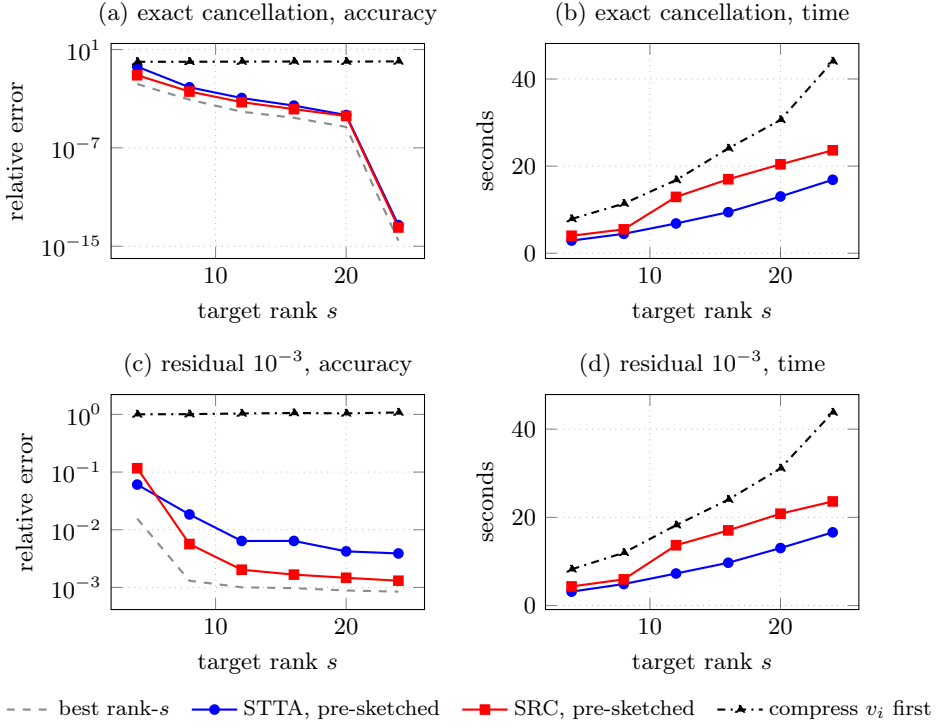

The two obvious ways of evaluating the sum are ruled out, each for its own
reason.  Accumulating the exact products fails on cost: a single $Av_i$ already
occupies $242$~MB, and a straightforward TT representation of their
accumulation would require about $2.4$~TB.  Keeping memory bounded by truncating each realization to rank $s$ and the running sum after every update, however, compromises accuracy. Since the stochastic contribution can dominate individual realizations, low-rank truncation may discard components of the deterministic field $u$ that subsequent averaging cannot recover.  As \Cref{fig:stream} shows, this route never achieves a relative
error below $1$ for any target rank tested, performing no better than the zero
approximation, and it is also the slowest of the routes considered.
 
What remains is to compress the product without ever forming it, and only once
the ensemble has been accumulated. This is what the pre-sketched STTA primitive
does: the sketches $\Omega_k$ and $\Psi_k$ are linear in the argument, so each
realization contributes to two accumulators and can then be discarded, and a
single recovery is performed at the end of the stream.  It remains within a
factor of $9$ to $26$ of the best rank-$s$ approximation of $Au$,
and is between $2$ and $3$ times faster than compressing the realizations
first.
 
Perfect cancellation of the stochastic perturbations is not to be expected in practice. For independent zero-mean perturbations with finite variance, their average decays like $N^{-1/2}$, whereas the deterministic part of the mean remains $u$. Thus, the stochastic contribution becomes progressively smaller as the ensemble grows, but does not vanish exactly. To assess the effect, we impose a residual stochastic contribution of relative size $10^{-3}$ in the ensemble mean. The solution then has TT rank $r^A r^v=224$ rather than $24$, and the error of its best rank-$s$ approximation plateaus at about $10^{-3}$ for $s\geq 8$. The pre-sketched approach follows this plateau, whereas compressing the realizations first retains essentially the same truncation bias as before.

For completeness, we also include pre-sketched SRC in the experiment. By exploiting the linearity of its first sweep, the first-sweep sketches of the realizations can be summed and orthogonalized once per site. Pre-sketched SRC is usually more accurate than pre-sketched STTA at equal target rank. However, its recovery requires access to the realizations again, either through multiple passes over the stream or by storing the entire ensemble. This entails a substantially higher memory cost, growing linearly with $N$, and can make pre-sketched SRC impractical when the ensemble does not fit in memory. Thus, when the realizations cannot be stored or revisited, pre-sketched STTA is the only suitable approach among those considered.
\section{Conclusions}
\label{sec:conclusions}

We have presented three randomized, streamable strategies for computing a
low-rank TT approximation of $Av$ without forming the product, based on
carefully ordered contractions and structured sketches within STTA. The
main contribution is the pre-sketched strategy, which absorbs the sketch
into the operator cores and is particularly effective for large mode sizes
when the target rank $s$ is smaller than the TT ranks of $v$. Our numerical
experiments on near-resonant linear systems and streaming uncertainty
quantification demonstrate the practical relevance of this regime, in which
the proposed method offers substantial computational advantages over
existing approaches with near-optimal accuracy guarantees.

We have also introduced a pre-sketched variant of SRC, the computational
cost of which is comparable to that of pre-sketched STTA but for which the use of orthogonal
rather than oblique projectors can yield greater accuracy. Among the methods
considered, it offers the most accurate alternative in this regime, but
requires access to the data beyond a single pass. Pre-sketched STTA therefore
provides an attractive combination of efficiency, accuracy, and streamability,
and can be the only viable option when the realizations cannot be stored or
revisited. For larger target ranks, standard SRC remains a competitive
alternative.

An interesting direction is to extend the pre-sketched approach to sparse
operators. While sparsity can already make the application of the operator
inexpensive, the orthogonalization steps in competing methods may limit their
computational advantage. Sketches that preserve sparsity in the operator
cores after sketch absorption could allow STTA to exploit this structure
throughout the computation and potentially outperform these methods.
Another open question is how aggressive compression can be combined with
operator recycling in iterative methods, as discussed in \cref{sec:krylov}.

\section*{Acknowledgments}
The authors acknowledge funding from the UK’s Engineering and Physical Sciences Research Council (EPSRC grant EP/Z533786/1). 
AI tools were used to improve the coding work, which the authors then independently checked and tested. The authors assume responsibility for all content.

\newpage

\bibliographystyle{siamplain}
\bibliography{references}
\end{document}